\documentclass [12pt]{article}
\usepackage[utf8]{inputenc}
\usepackage[english]{babel}
\usepackage{amsmath}
\usepackage{amsfonts, amssymb}
\usepackage[left=2.5cm,right=2.2cm,top=2.5cm,bottom=3.5cm]{geometry}
\usepackage{comment}
\usepackage{cmap}

\let\originalleft\left
\let\originalright\right
\renewcommand{\left}{\mathopen{} \mathclose\bgroup\originalleft}
\renewcommand{\right}{\aftergroup\egroup\originalright}

\newenvironment{proof}[1][] {\noindent {\bf Proof#1:} }{\hspace*{\fill}$\square$\medskip\par}

\newtheorem{thm}{Theorem}[section]
\newtheorem{lem}[thm]{Lemma}
\newtheorem{prop}[thm]{Proposition}
\newtheorem{cor}[thm]{Corollary}

\newtheorem{rem}[thm]{Remark}

\newtheorem{cond}[thm]{Condition}

\author{S.E.Nikitin\thanks{St. Petersburg Department of Steklov Mathematical Institute of Russian Academy of Sciences, 27 Fontanka, St. Petersburg, 191023, Russia. Email: \texttt{nikitin97156@mail.ru}.}}
\title{Boundary gradient asymptotics for superquadratic viscous Hamilton--Jacobi equations via optimal control\thanks{This work was supported by the Ministry of Science and Higher Education of the Russian Federation (agreement 075-15-2025-344 dated 29/04/2025 for Saint Petersburg Leonhard Euler International Mathematical Institute at PDMI RAS).}}	
\date{November 2026}
\medskip

\begin{document}
\maketitle

\begin{abstract}
We study the boundary behavior of the maximal solution of
\[
	-\Delta u + |\nabla u|^{p} + \lambda u = f
\]
on arbitrary bounded $C^{2}$ domains, where $p > 2$ and $\lambda > 0$. For bounded state-dependent terms $f$, we resolve a conjecture of Lasry and Lions by establishing sharp boundary gradient asymptotics without geometric restrictions on the domain. Previously, such asymptotics were established only for certain classes of domains with additional restrictions on $p$ and $f$. We also extend the analysis to state-dependent terms that may be singular near the boundary and obtain sharp boundary profiles and gradient asymptotics.

Our approach relies on stochastic optimal control and boundary blow-up analysis. Dynamic programming estimates and comparisons of near-optimal trajectories establish tangential invariance of boundary blow-up limits, leading to the identification of one-dimensional profiles and the asymptotic formulas. For the associated state-constrained control problem, we establish sharp asymptotics for the minimizing feedback and deduce its admissibility and optimality, with corresponding results for the ergodic problem.
\end{abstract}



\section{Introduction}

Optimal control problems with state constraints lead naturally to Hamilton--Jacobi--Bellman equations with singular boundary behavior. A basic model studied systematically by Lasry and Lions \cite{LL89} is
\begin{equation} \label{e:u}
- \Delta u + \lvert \nabla u \rvert^{p} + \lambda u = f
\qquad \text{in } \Omega,
\end{equation}
where $\Omega \subset \mathbb{R}^{N}$ is a bounded $C^{2}$ domain, $p > 1$, $\lambda > 0$, and $f$ is a given function. In the associated stochastic control problem,
\[
\mathrm{d}X_{t} = b_{t} \, \mathrm{d}t + \sqrt{2} \, \mathrm{d}W_{t}, \qquad t \ge 0,
\]
the admissible controls are required to keep $X_{t}$ in $\Omega$ for all $t \ge 0$. The running cost is $(p - 1) p^{-q} \lvert b_{t} \rvert^{q} + f(X_{t})$, where $q := p / (p - 1)$, and the formal minimizing feedback is
\begin{equation} \label{def:b0}
b^{*}(x) := -p \lvert \nabla u(x) \rvert^{p - 2} \nabla u(x), \qquad x \in \Omega.
\end{equation}

For bounded $f$, the maximal solution of \eqref{e:u} constructed by Lasry and Lions blows up at the boundary when $1 < p \le 2$, whereas it remains continuous up to the boundary for $p > 2$ despite the singular behavior of its gradient. In this paper, we consider only the superquadratic regime $p > 2$. The sharp boundary gradient asymptotics established by Lasry and Lions under additional assumptions take the form
\[
d(x)^{q - 1} \nabla u(x) \longrightarrow (q - 1)^{q - 1} n(\pi_{x}) \qquad \text{as } d(x) \downarrow 0,
\]
where $d(x) = \operatorname{dist}(x, \partial \Omega)$, $\pi_{x}$ is the nearest point of $\partial \Omega$ to $x$, and $n$ is the outward unit normal. Consequently,
\[
d(x) b^{*}(x) \longrightarrow -q \, n(\pi_{x}).
\]
The resulting inward drift has magnitude asymptotic to $q d^{-1}$ and, since $q > 1$, prevents the controlled diffusion from reaching the boundary.

Lasry and Lions established the sharp boundary asymptotics for $f \in W^{1, \infty}(\Omega)$ under additional geometric assumptions. Their Proposition~IV.3 covers balls and half-spaces for all $p > 2$, and convex domains for $2 < p < 3$. They also observed in Remark~IV.6 that, for convex domains, the same argument extends to all $p > 2$ when $f$ itself is convex.

Their boundary blow-up argument relies on tangential regularity estimates to show that limiting profiles depend only on the normal variable, a property we call {\it tangential invariance}. These estimates are not available in general domains, and Lasry and Lions conjectured \cite[p.~608]{LL89} that the sharp asymptotics hold for bounded $f$ without geometric restrictions. The first purpose of this paper is to prove this conjecture on arbitrary bounded $C^{2}$ domains.

Restrictions of the form $p \le 3$ also appear in related time-dependent superquadratic viscous Hamilton--Jacobi equations with Dirichlet boundary conditions. Porretta and Souplet \cite{PS17} established, for $2 < p \le 3$, sharp tangential gradient blow-up profiles in two space dimensions for solutions exhibiting single-point gradient blow-up on a locally flat portion of the boundary, and left the case $p > 3$ open. Similar restrictions arise in \cite{AS20, CHZ23}. In contrast, Filippucci, Pucci, and Souplet \cite{FPS20} established sharp normal gradient blow-up profiles for the time-dependent problem and optimal boundary gradient estimates for the elliptic problem, both with homogeneous Dirichlet conditions and for all $p > 2$ in general smooth domains. Their half-space Liouville theorem relies on homogeneous Dirichlet conditions, whereas in our setting the boundary values of the blow-up limits are not known a priori.

We also go beyond the bounded-data conjecture by allowing the right-hand side to blow up near the boundary, with $f = a d^{-\beta} + o(d^{-\beta})$ for $q \le \beta < p$. For $\beta = q$, the case $a = 0$ on $\partial \Omega$ is {\it subcritical} and directly extends the bounded-data setting of Lasry and Lions, while positive values of $a$ give rise to the {\it critical} regime. In the subcritical case, lower-order terms disappear under the natural boundary rescaling, whereas in the critical case a term of order $d^{-q}$ survives and changes the limiting equation. We establish the corresponding sharp gradient asymptotics and allow the leading coefficient $a$ to vary along the boundary.

In the {\it supercritical} regime $q < \beta < p$, the diffusion term disappears under the natural rescaling, and the limiting profile solves a first-order eikonal equation. Assuming $a > 0$ on $\partial \Omega$, we prove
\[
d(x)^{\beta / p} \nabla u(x) \longrightarrow a(\pi_{x})^{1 / p} n(\pi_{x}), \qquad d(x)^{\beta / q} b^{*}(x) \longrightarrow -p \, a(\pi_{x})^{1 / q} n(\pi_{x}),
\]
so the inward drift has magnitude of order $d^{-\beta / q}$, which grows faster than $d^{-1}$. The corresponding boundary profile asymptotics in all three regimes were anticipated by Lasry and Lions \cite[Rem.~IV.5]{LL89}. For $\beta \ge p$, solutions bounded from below blow up at the boundary; this regime is treated in \cite[Sect.~III]{LL89} and is not considered here.

Our approach uses the Lasry--Lions identification of the maximal solution with the value function of the associated state-constrained control problem. Dynamic programming inequalities and comparisons of near-optimal controlled trajectories establish tangential invariance of boundary blow-up limits without geometric assumptions or tangential regularity estimates. The sharp gradient asymptotics also yield admissibility and optimality of the minimizing feedback.

Dutta, Nguyen, and Tu \cite{DNT26} recently studied the vanishing-viscosity limit of maximal solutions of superquadratic Hamilton--Jacobi equations with state constraints, using boundary gradient estimates as an auxiliary ingredient. Here we prove the corresponding sharp boundary gradient asymptotics. We also establish a general local $C^{1}$-convergence result for vanishing-viscosity approximations, of independent interest.

The paper is organized as follows. Section \ref{s:smr} states the main results. Section \ref{s:p} develops the stochastic-control preliminaries, dynamic programming inequalities, and boundary blow-up limits. Section \ref{s:ti} establishes estimates for near-optimal controls and proves tangential invariance. Section \ref{s:lo} identifies the limiting one-dimensional profiles in the critical and supercritical cases and proves local $C^{1}$ convergence of the rescaled solutions. Finally, Section \ref{s:cop} completes the proofs of the main results.

\section{Setting and main results} \label{s:smr}
Let $\Omega \subset \mathbb{R}^{N}$ be a bounded $C^{2}$ domain with outward unit normal $n$. For $\delta > 0$, let $\Omega_{\delta} := \{x \in \Omega : \operatorname{dist}(x, \partial \Omega) < \delta\}$. Using \cite[Lem.~14.16]{GT01}, we fix $\delta_{0} \in (0, 1)$ such that $\operatorname{dist}(\cdot, \partial \Omega)$ is $C^{2}$ on $\overline{\Omega}_{2 \delta_{0}}$ and the nearest-point projection $\pi : \Omega_{2 \delta_{0}} \to \partial \Omega$ is well-defined. We write $\pi_{x} := \pi(x)$. We fix $d \in C^{2}(\overline{\Omega})$ with $d = \operatorname{dist}(\cdot, \partial \Omega)$ on $\Omega_{\delta_{0}}$ and $d \ge \delta_{0}$ on $\Omega \setminus \Omega_{\delta_{0}}$. Then $\nabla d(x) = - n(\pi_{x})$ for $x \in \Omega_{\delta_{0}}$. 

Let
\[
p > 2, \qquad q := \frac{p}{p - 1} \in (1, 2), \qquad \lambda > 0, \qquad \beta \in [q, p)
\]
and set
\[
\alpha := 1 - \frac{\beta}{p}, \qquad \theta := \frac{\beta}{p} = 1 - \alpha.
\]
For $\beta = q$, these are $\alpha = 2 - q = \frac{p - 2}{p - 1}$ and $\theta = q - 1 = \frac{1}{p - 1}$.

Let $a(x) := d(x)^{\beta} f(x)$ on $\Omega_{\delta_{0}}$. We consider the following two cases, distinguished by the boundary growth rate of $f$.
\begin{cond}[Subcritical and critical case] \label{cond:critical}
$f \in L^{\infty}_{\mathrm{loc}}(\Omega)$;
$f \ge -c_{f}$ a.e. in $\Omega$ for some $c_{f} \ge 0$;
$\beta = q$, and there exists $\bar{a} \in C(\partial\Omega)$ such that
\begin{equation} \label{e:critical-trace}
\lim_{\delta \downarrow 0}
\left\| a - \bar{a} \circ \pi \right\|_{L^{\infty}(\Omega_\delta)} = 0.
\end{equation}
\end{cond}
\begin{cond}[Supercritical case] \label{cond:supercritical}
$f \in W^{1, \infty}_{\mathrm{loc}}(\Omega)$;
$f \ge -c_{f}$ a.e. in $\Omega$ for some $c_{f} \ge 0$; $\beta \in (q, p)$ and $a(\cdot)$ extends continuously to $\overline{\Omega}_{\delta_{0}}$ with $a > 0$ on $\partial \Omega$; and
\begin{equation} \label{e:wdf}
\operatorname*{ess\,sup}_{\Omega_{\delta_{0}}}
 d^{\beta + 1} \lvert \nabla f \rvert < \infty.
\end{equation}
\end{cond}
We always choose a Borel representative of $f$ such that $f \ge -c_{f}$ everywhere in $\Omega$, and under Condition~\ref{cond:supercritical} we take the continuous representative of $f$. By the occupation estimate of \cite[Th.~4.1]{WLW23}, modifying $f$ on a set of Lebesgue measure zero does not affect the time integrals of $f$ along controlled trajectories.
Under Condition~\ref{cond:critical}, we set $a := \bar{a}$ on $\partial \Omega$. Thus, under either condition, $a$ has continuous boundary values, which are nonnegative due to $f \ge -c_{f}$. Moreover, for some $C_f > 0$,
\begin{equation} \label{e:cg}
\lvert f(x) \rvert \le C_{f} \bigl(1 + d(x)^{-\beta} \bigr), \qquad \text{for a.e. } x \in \Omega.
\end{equation}
We call Condition~\ref{cond:critical} with $a \equiv 0$ on $\partial \Omega$ subcritical. In particular, this includes every $f$ bounded from below such that $f^{+}(x) = O(d(x)^{-\beta^{\prime} })$ as $d(x) \to 0$ with $\beta^{\prime} \in [0,q)$.

Let $u$ denote the maximal solution of \eqref{e:u}, that is, the largest solution in $W^{2, k}_{\mathrm{loc}}(\Omega)$ for all $k \in [1, \infty)$. By \cite[Cor.~IV.1 and IV.2]{LL89}, such a solution exists and belongs to $C(\overline{\Omega})$.
We say that $v \in C(\overline{\Omega})$ satisfies the state-constraint boundary condition if, for every $\psi \in C^{2}(\overline{\Omega})$,
\begin{equation} \label{e:state-constraint}
\text{$v - \psi$ has no local minimum on $\overline{\Omega}$ at a point of $\partial \Omega$.}
\end{equation}
By \cite[Sect.~V.1 and (62)]{LL89}, the maximal solution $u$ satisfies \eqref{e:state-constraint}.

Our main results give both the boundary profile of $u$ and the corresponding sharp gradient asymptotics. For $s \ge 0$, let $A(s)$ be the unique root in
$[(q - 1)^{q - 1}, \infty)$ of
\begin{equation} \label{e:critical-algebraic-B}
A^{p} - (q - 1)A = s.
\end{equation}
The function $A$ is continuous and increasing on $[0, \infty)$, with $A(0) = (q - 1)^{q - 1}$.

\begin{thm}
\label{t:c}
Under Condition~\ref{cond:critical}, $u$ satisfies
\begin{equation} \label{e:critical-main-profile}
\lim_{t \downarrow 0} \sup_{x_{0} \in \partial \Omega} \left|
 t^{-\alpha} \bigl(u(x_{0} - t n(x_{0})) - u(x_{0}) \bigr) + \frac{A(a(x_{0}))}{\alpha} \right| = 0,
\end{equation}
and
\begin{equation} \label{e:cmg}
\lim_{\varepsilon \downarrow 0} \sup_{\{0 < d(x) < \varepsilon\}} \left|
 d(x)^{\theta} \nabla u(x) - A(a(\pi_{x})) n(\pi_{x}) \right| = 0.
\end{equation}
\end{thm}
\begin{cor} \label{c:main}
Assume Condition~\ref{cond:critical} with $a \equiv 0$ on $\partial \Omega$. Then
\[
d(x)^{q - 1} \nabla u(x) \longrightarrow (q - 1)^{q - 1} n(\pi_{x}) \qquad \text{uniformly as } d(x) \downarrow 0.
\]
The corresponding boundary profile is given by \eqref{e:critical-main-profile} with $A(0) = (q - 1)^{q - 1}$. In particular, the result holds for $f \in L^{\infty}(\Omega)$, which is the conjecture of \cite[p.~608]{LL89}.
\end{cor}
\begin{thm}
\label{t:supc}
Under Condition~\ref{cond:supercritical}, $u$ satisfies
\begin{equation*} 
\lim_{t \downarrow 0} \sup_{x_{0} \in \partial \Omega} \left|
 t^{-\alpha} \bigl(u(x_{0} - t n(x_{0})) - u(x_{0}) \bigr) + \frac{a(x_{0})^{1 / p}}{\alpha} \right| = 0,
\end{equation*}
and
\begin{equation} \label{e:smg}
\lim_{\varepsilon \downarrow 0} \sup_{\{0 < d(x) < \varepsilon\}} \left|
 d(x)^{\theta} \nabla u(x) - a(\pi_{x})^{1 / p} n(\pi_{x}) \right| = 0.
\end{equation}
\end{thm}
\begin{rem}
The boundary asymptotics in Theorems~\ref{t:c} and \ref{t:supc} are local in nature. The arguments below can be localized to establish the corresponding asymptotic formulas locally uniformly along a relatively open $C^{2}$ portion $\Gamma$ of $\partial \Omega$ near which Condition~\ref{cond:critical} or \ref{cond:supercritical} holds, provided that a bounded maximal solution exists and admits a control representation (see Subsection~\ref{ss:control} and Proposition~\ref{p:verification}). In particular, different growth regimes may occur on different portions of the boundary, and $\Omega$ need not be bounded. We do not provide the precise conditions or the details of this localization here.
\end{rem}

\subsection{Results for control problems} \label{ss:control}
We next state the consequences for the stochastic control problem. Set
\begin{equation*}
L(b) := (p - 1) p^{-q} \lvert b \rvert^{q}, \qquad b \in \mathbb{R}^{N}.
\end{equation*}
For an $N$-dimensional Brownian motion $W$ and a progressively measurable control process $b$ such that $\int_{0}^{t} \lvert b_{s} \rvert \, \mathrm{d}s < \infty$ a.s. for every $t < \infty$, consider
\begin{equation} \label{e:X} 
\mathrm{d}X_{t} = b_{t} \, \mathrm{d}t + \sqrt{2} \, \mathrm{d}W_{t},
\qquad X_{0} = x \in \Omega.
\end{equation}

We work in the weak formulation: the filtered probability space, the Brownian motion $W$, the control $b$ and the corresponding process $X$ are part of the control system. Such a system is admissible for $x \in \Omega$ if \eqref{e:X} holds and, almost surely, $X_{t} \in \Omega$ for every $t \ge 0$. We denote the class of admissible control systems by $\operatorname{Adm}(x)$ and write $b \in \operatorname{Adm}(x)$ for brevity. For admissible controls, define the cost functional and the value function 
\begin{equation*} 
J_{\lambda}(x; b) := \mathbb{E} \int_{0}^{\infty} e^{-\lambda t}
\bigl(L(b_{t}) + f(X_{t}) \bigr) \, \mathrm{d}t, \qquad 
V_{\lambda}(x) := \inf_{b \in \operatorname{Adm}(x)} J_{\lambda}(x; b).
\end{equation*}
Since $L \ge 0$ and $f \ge -c_{f}$, $J_{\lambda}(x; b)$ is well defined with values in $[-c_{f} / \lambda, \infty]$ for every $x \in \Omega$ and $b \in \operatorname{Adm}(x)$. The relation between this control problem and \eqref{e:u} is based on the convex duality identity
\begin{equation} \label{e:Fenchel-duality}
\inf_{b \in \mathbb{R}^{N}} \{ L(b) + b \cdot v \} = - \lvert v \rvert^{p}, \qquad v \in \mathbb{R}^{N}.
\end{equation}

The sharp gradient asymptotics ensure the admissibility of the feedback 
$b^{*}$ defined by \eqref{def:b0}. Set
\begin{equation*}
\mu(s) := p A(s)^{p - 1} \ge q, \qquad s \ge 0.
\end{equation*}
\begin{thm} \label{t:critical-feedback}
Under Condition~\ref{cond:critical},
$b^{*} \in C(\Omega; \mathbb{R}^{N})$ and 
\[
d(x) b^{*}(x) \longrightarrow - \mu(a(\pi_{x})) \, n(\pi_{x}) \qquad \text{as } d(x) \to 0.
\] 
For every $x \in \Omega$ there exists a weak solution, unique in law, of
\begin{equation} \label{e:feedback-SDE}
\mathrm{d}X_{t} = b^{*}(X_{t}) \, \mathrm{d}t + \sqrt{2} \, \mathrm{d}W_{t},
 \qquad X_{0} = x,
\end{equation}
and it remains in $\Omega$ for all $t \ge 0$ almost surely. Moreover,
\begin{equation} \label{e:critical-feedback-optimality}
J_{\lambda}(x; b^{*}(X)) = V_{\lambda}(x) = u(x).
\end{equation}
Thus $b^{*}$ is an optimal Markov feedback. 
\end{thm}

\begin{thm}
\label{t:supercritical-feedback}
Under Condition~\ref{cond:supercritical}, $b^{*} \in C(\Omega; \mathbb{R}^{N})$ and
\[
d(x)^{\beta / q} b^{*}(x) \longrightarrow - p \, a(\pi_{x})^{1 / q} n(\pi_{x}) \qquad \text{as } d(x) \to 0.
\]
Moreover, the remaining conclusions of Theorem~\ref{t:critical-feedback} hold: for every $x \in \Omega$, \eqref{e:feedback-SDE} has a weak solution, unique in law, which remains in $\Omega$ for all $t \ge 0$ almost surely, and \eqref{e:critical-feedback-optimality} holds.
\end{thm}
In particular, under either condition, the minimizing feedback is admissible and optimal on every bounded $C^{2}$ domain, without the geometric assumptions of \cite[Prop.~IV.3]{LL89}.

\subsection{The ergodic problem}
We also consider the ergodic equation
\begin{equation} \label{e:ergodic}
-\Delta u_{0} + \lvert \nabla u_{0} \rvert^{p} + \kappa = f \qquad \text{in } \Omega,
\end{equation}
where $\kappa \in \mathbb{R}$ is an unknown constant. By \cite[Th.~VI.3]{LL89}, under Condition~\ref{cond:critical} or \ref{cond:supercritical}, there exists a unique ergodic constant $\kappa$ and a solution $u_{0} \in W^{2, k}_{\mathrm{loc}}(\Omega) \cap C(\overline{\Omega})$ for every $k \in [1, \infty)$, satisfying \eqref{e:state-constraint} and unique up to an additive constant among solutions satisfying this boundary condition. More precisely, fix any $\bar{x} \in \Omega$ and normalize $u_{0}$ by $u_{0}(\bar{x}) = 0$. If $u_{\lambda}$ denotes the maximal solution of \eqref{e:u} for the discount factor $\lambda > 0$, then
\[
\lambda u_{\lambda} \longrightarrow \kappa, \qquad u_{\lambda} - u_{\lambda}(\bar{x}) \longrightarrow u_{0} \qquad \text{uniformly on } \overline{\Omega} \text{ as } \lambda \downarrow 0.
\]

Our boundary gradient asymptotics and feedback optimality results extend to the ergodic problem.

\begin{cor} \label{c:ergodic}
Under Condition~\ref{cond:critical} or \ref{cond:supercritical}, the ergodic solution $u_{0}$ satisfies the conclusions of Theorems~\ref{t:c} and \ref{t:supc}, respectively, with $u_{0}$ in place of $u$.

Moreover, the feedback
\[
b_{0}^{*}(x) := -p \lvert \nabla u_{0}(x) \rvert^{p - 2} \nabla u_{0}(x), \qquad x \in \Omega
\]
is admissible and optimal for the long-run average cost. More precisely, for every $x \in \Omega$,
\begin{equation} \label{e:ergodic-value}
\kappa = \inf_{b \in \operatorname{Adm}(x)} \liminf_{T \to \infty} \frac{1}{T} \mathbb{E} \int_{0}^{T} \bigl(L(b_{t}) + f(X_{t})\bigr) \, \mathrm{d}t,
\end{equation}
and the infimum is attained by $b_{t} = b_{0}^{*}(X_{t})$, where $X$ satisfies \eqref{e:feedback-SDE} with $b_{0}^{*}$ in place of $b^{*}$. 
\end{cor}

\section{Preliminaries and boundary blow-up} \label{s:p}

From now on, we assume Condition~\ref{cond:critical} or \ref{cond:supercritical}. 

\subsection{Lasry--Lions inputs}
\begin{prop} \label{p:LL}
For all $k \in [1, \infty)$,
\begin{equation} \label{e:LL-regularity}
u \in W^{2, k}_{\mathrm{loc}}(\Omega) \cap C(\overline{\Omega}).
\end{equation}
Moreover, for some $C_{H} > 0$,
\begin{equation} \label{e:LL-holder}
\lvert u(x) - u(y) \rvert \le C_{H} \lvert x - y \rvert^{\alpha}, \qquad x, y \in \overline{\Omega},
\end{equation}
and, for some $C_{G} > 0$,
\begin{equation} \label{e:LL-gradient}
\lvert \nabla u(x) \rvert \le C_{G} \, d(x)^{-\theta}, \qquad x \in \Omega.
\end{equation}
Finally, for every $x \in \Omega$ and $\varepsilon > 0$, there exists $b^{\varepsilon} \in \operatorname{Adm}(x)$ with
\begin{equation} \label{e:LL-value}
J_{\lambda}(x; b^{\varepsilon}) < u(x) + \varepsilon.
\end{equation}
\end{prop}
\begin{proof}
Under Condition~\ref{cond:critical} or \ref{cond:supercritical}, the growth bound \eqref{e:cg} and, in the supercritical case, the derivative bound \eqref{e:wdf} ensure that the results of \cite[Th.~IV.1, Cor.~IV.1, IV.2]{LL89} apply, yielding \eqref{e:LL-regularity}--\eqref{e:LL-gradient}.
Finally, by \cite[Sect.~V.1, Th.~VII.2]{LL89}, $u$ coincides with the value function over admissible Markov feedbacks, which yields \eqref{e:LL-value}.
\end{proof}

\subsection{Dynamic programming inequalities}
\begin{prop}
\label{p:verification}
Let $x \in \Omega$, let $b$ be a progressively measurable process with $\int_{0}^{t} \lvert b_{s} \rvert \, \mathrm{d}s < \infty$ a.s. for every $t \ge 0$, and let $X$ be a corresponding solution to \eqref{e:X}. If $\tau$ is a bounded stopping time such that, almost surely, $X_{t} \in \Omega$ for all $t \in [0, \tau]$, then
\begin{equation} \label{e:uu}
u(x) \le \mathbb{E} \biggl(\int_{0}^{\tau} e^{-\lambda t} \bigl(L(b_{t}) + f(X_{t}) \bigr) \, \mathrm{d}t + e^{-\lambda \tau} u(X_{\tau}) \biggr).
\end{equation}
For every $\varepsilon > 0$, there exists $b^{\varepsilon} \in \operatorname{Adm}(x)$, with corresponding controlled process $X^{\varepsilon}$, such that, for every bounded stopping time $\sigma$ of the filtration of the control system,
\begin{equation} \label{e:dppr}
\mathbb{E} \biggl(\int_{0}^{\sigma} e^{-\lambda t} \bigl(L(b_{t}^{\varepsilon}) + f(X_{t}^{\varepsilon}) \bigr) \, \mathrm{d}t + e^{-\lambda \sigma} u(X_{\sigma}^{\varepsilon}) \biggr) \le J_{\lambda}(x; b^{\varepsilon}) < u(x) + \varepsilon.
\end{equation}
In particular, $u = V_{\lambda}$ in $\Omega$.
\end{prop}

\begin{proof}
By \eqref{e:Fenchel-duality}, we have
$L(b) + b \cdot v + \lvert v \rvert^{p} \ge 0$ for $b, v \in \mathbb{R}^{N}$. Therefore, a.e. in $\Omega$, 
\begin{equation} \label{e:generator-inequality}
(\Delta + b \cdot \nabla - \lambda) u \ge - L(b) - f, \qquad b \in \mathbb{R}^{N}.
\end{equation}
Let $\sigma$ be a stopping time such that $\sigma \le \tau$. We will first prove that
\begin{equation} \label{e:st}
\mathbb{E} \, e^{-\lambda \sigma} u(X_{\sigma}) \le \mathbb{E} \biggl(
\int_{\sigma}^{\tau} e^{-\lambda t} \bigl(L(b_{t}) + f(X_{t}) \bigr) \, \mathrm{d}t
+ e^{-\lambda \tau} u(X_{\tau}) \biggr).
\end{equation}

For $i, j \in \mathbb{N}$, set $\tau_{i} := \tau \, \wedge \, \inf \{t \ge \sigma : d(X_{t}) \le i^{-1} \}$ and $\nu_{j} := \inf \left\{t \ge 0 : \int_{0}^{t} \lvert b_{s} \rvert \, \mathrm{d}s \ge j \right\}$. By \eqref{e:LL-regularity}, $u \in W^{2, k}_{\mathrm{loc}}(\Omega)$ for all $k \in [1, \infty)$, so the It\^{o}--Krylov formula applies on $[\sigma \wedge \nu_{j}, \tau_{i} \wedge \nu_{j}]$. Indeed, it follows by smooth approximation of $u$ and the occupation estimate of \cite[Th.~4.1]{WLW23} applied up to $\tau_{i} \wedge \nu_{j}$, which also justifies applying the almost-everywhere inequality \eqref{e:generator-inequality} along $X$. The localization ensures that the resulting stochastic integral is a true martingale and hence has zero expectation. Taking expectations and using \eqref{e:generator-inequality} gives \eqref{e:st} on $[\sigma \wedge \nu_{j}, \tau_{i} \wedge \nu_{j}]$ instead of $[\sigma, \tau]$. Letting first $i \to \infty$ and then $j \to \infty$, monotone convergence for $L(b) + f + c_{f} \ge 0$ and dominated convergence for bounded endpoint terms give \eqref{e:st}, because, almost surely, $\tau_{i} = \tau$ for all sufficiently large $i$ and $\nu_{j} \uparrow \infty$ as $j \to \infty$. With $\sigma = 0$, \eqref{e:st} is \eqref{e:uu}.

For $b \in \operatorname{Adm}(x)$, applying \eqref{e:uu} with $\tau = T$ and letting $T \to \infty$, we obtain $u(x) \le J_{\lambda}(x; b)$, since $u$ is bounded. Thus $u \le V_{\lambda}$. By \eqref{e:LL-value},
\begin{equation*}
u = V_{\lambda} \qquad \text{in } \Omega. 
\end{equation*}
For $\varepsilon > 0$, let $b^{\varepsilon} \in \operatorname{Adm}(x)$ be such that $J_{\lambda}(x; b^{\varepsilon}) < u(x) + \varepsilon$. Let $\sigma$ be as in \eqref{e:dppr}.
Applying \eqref{e:st} to $b^{\varepsilon}$, $X^{\varepsilon}$ with $\tau = T$ for a deterministic constant $T \ge \sigma$, adding the running cost on $[0, \sigma]$, and letting $T \to \infty$, we obtain 
\[
\mathbb{E} \biggl(\int_{0}^{\sigma} e^{-\lambda t} \bigl(L(b_{t}^{\varepsilon}) + f(X_{t}^{\varepsilon}) \bigr) \, \mathrm{d}t + e^{-\lambda \sigma} u(X_{\sigma}^{\varepsilon}) \biggr) \le J_{\lambda}(x; b^{\varepsilon}) < u(x) + \varepsilon,
\]
which is \eqref{e:dppr}.
\end{proof}

\subsection{Scaling and blow-up limits}
Fix $x_{0} \in \partial \Omega$ and $r_{j} \downarrow 0$. After a translation and a rotation, under which we keep the notation $u$, $f$, $a$ and $d$ for the transformed functions, assume $x_{0} = 0$ and the inward unit normal at $0$ is $e_{N}$. Write $\xi = (z, \rho) \in \mathbb{R}^{N - 1} \times \mathbb{R}$. For a sufficiently small cylinder $\mathcal{C}_{0}$ centered at $0$, the domain is represented locally as
\[
\Omega \cap \mathcal{C}_{0} = \{(z, \rho) \in \mathcal{C}_{0} : \rho > g(z) \},
\]
for some $C^{2}$ function $g$ defined on the base of $\mathcal{C}_{0}$ and satisfying, for some $C_{g} > 0$,
\begin{equation} \label{e:g}
 g(0) = 0, \qquad \nabla g(0) = 0, \qquad \lvert D^{2} g \rvert \le C_{g}.
\end{equation}
Define the rescaled domain and the rescaled solution as
\begin{equation*}
\Omega_{j} := \frac{\Omega}{r_{j}}, \qquad u_{j}(\xi) := r_{j}^{-\alpha} \bigl(u(r_{j} \xi) - u(0) \bigr), \quad \xi \in \overline{\Omega}_{j},
\end{equation*}
and set $g_{j}(z) := g(r_{j} z) / r_{j}$ and $H := \{(z, \rho) : \rho > 0\}$.
The estimate \eqref{e:LL-holder} is scale invariant:
\begin{equation} \label{e:scaled-holder}
\lvert u_{j}(\xi) - u_{j}(\eta) \rvert \le C_{H} \lvert \xi - \eta \rvert^{\alpha}, \qquad \xi, \eta \in \overline{\Omega}_{j}.
\end{equation}
Since $u_{j}(0) = 0$, we have
\begin{equation} \label{e:scaled-growth}
\lvert u_{j}(\xi) \rvert \le C_{H} \lvert \xi \rvert^{\alpha}, \qquad \xi \in \overline{\Omega}_{j}.
\end{equation}
The next lemma establishes the convergence of the rescaled domains to the half-space and the existence of blow-up limits.
\begin{lem} \label{l:blowup}
For every bounded cylinder $\mathcal{C}$ centered at $0$ and for all large $j$,
\[
\Omega_{j} \cap \mathcal{C} = \{(z, \rho) \in \mathcal{C} : \rho > g_{j}(z)\},
\]
where, for $z$ in the base of $\mathcal{C}$,
\begin{equation} \label{e:scaled-graph}
\lvert \nabla g_{j}(z) \rvert \le C_{g} r_{j} \lvert z \rvert, \qquad \lvert g_{j}(z) \rvert \le \frac{C_{g}}{2} r_{j} \lvert z \rvert^{2}.
\end{equation}
In particular, $\Omega_{j} \to H$ locally: every compact subset of $H$ is contained in $\Omega_{j}$ for all large $j$.
Moreover, after passing to a subsequence,
\begin{equation} \label{e:blowup-limit}
u_{j} \longrightarrow \Phi \qquad \text{locally uniformly on } H,
\end{equation}
where $\Phi$ extends uniquely to $\overline{H}$ and satisfies
\[
\lvert \Phi(\xi) - \Phi(\eta) \rvert \le C_{H} \lvert \xi - \eta \rvert^{\alpha}, \qquad \lvert \Phi(\xi) \rvert \le C_{H} \lvert \xi \rvert^{\alpha}, \qquad \xi, \eta \in \overline{H}.
\]
We call such $\Phi$ a blow-up limit at $x_{0}$.
\end{lem}

\begin{proof}
For every bounded cylinder $\mathcal{C}$ centered at $0$, we have $\Omega_{j} \cap \mathcal{C} = r_{j}^{-1} (\Omega \cap r_{j} \mathcal{C})$ and $r_{j} \mathcal{C} \subset \mathcal{C}_{0}$ for large $j$, which gives the graph representation. Since $\nabla g_{j}(z) = \nabla g(r_{j} z)$, the bounds \eqref{e:scaled-graph} follow from \eqref{e:g}. In particular, $g_{j} \to 0$ uniformly on the base of $\mathcal{C}$, so every compact subset of $H \cap \mathcal{C}$ lies in $\Omega_{j}$ for large $j$.

By \eqref{e:scaled-holder} and \eqref{e:scaled-growth}, the $u_{j}$ are equicontinuous and uniformly bounded on every compact subset of $H$, so Arzel\`{a}--Ascoli and a diagonal argument give \eqref{e:blowup-limit}. Both bounds pass to the limit, and the H\"{o}lder bound extends $\Phi$ uniquely to $\overline{H}$.
\end{proof}

\section{Tangential invariance} \label{s:ti}

\subsection{Scaled near-optimal controls: energy and localization} \label{ss:snoc}
Fix $\xi \in H$, set $x_{j} := r_{j} \xi$, and let $\varepsilon_j \downarrow 0$; then $x_{j} \in \Omega$ for all sufficiently large $j$. By Proposition~\ref{p:verification} applied at $x_{j}$ with $\varepsilon = r_{j}^{\alpha} \varepsilon_{j}$, choose an admissible control $b^{j}$, with controlled process $X^{j}$ and Brownian motion $W^{j}$, satisfying \eqref{e:dppr}. The control $b^{j}$ is chosen once for each $j$ and does not depend on the later parameters $T$ and $M$.
Set $m := \alpha + \beta$ and 
\[
Y_{s}^{j} := \frac{X_{r_{j}^{m} s}^{j}}{r_{j}}, \qquad 
\widetilde{b}_{s}^{j} := r_{j}^{m - 1} b_{r_{j}^{m} s}^{j}, \qquad \widetilde{W}_{s}^{j} := r_{j}^{-m / 2} W^{j}_{r_{j}^{m} s}, \qquad s \ge 0.
\]
Then, with $\varsigma_{j} := r_{j}^{m - 2}$,
\begin{equation*}
\mathrm{d}Y_{s}^{j} = \widetilde{b}_{s}^{j} \, \mathrm{d}s + \sqrt{2 \varsigma_{j}} \, \mathrm{d} \widetilde{W}_{s}^{j}, \qquad Y_{s}^{j} \in \Omega_{j}, \quad Y_{0}^{j} = \xi, \qquad s \ge 0.
\end{equation*}
Since $m = 1 + \beta / q$, we have $2 \le m < p$, with $m = 2$ iff $\beta = q$; thus $\varsigma_{j} = 1$ under Condition~\ref{cond:critical} and $\varsigma_{j} \to 0$ under Condition~\ref{cond:supercritical}.
Define the rescaled state cost 
\begin{equation} \label{e:scaled-forcing}
F_{j}(y) := r_{j}^{\beta} (f(r_{j} y) - \lambda u(0)), \qquad y \in \Omega_{j},
\end{equation}
and the rescaled distance function
\[
d_{j}(y) := \frac{d(r_{j} y)}{r_{j}}, \qquad y \in \Omega_{j}.
\]

For every bounded stopping time $\tau$ in the natural filtration of $Y^{j}$, applying \eqref{e:dppr} with $\sigma = r_{j}^{m} \tau$ gives
\begin{equation} \label{e:scaled-stopped}
u_{j}(\xi) + \varepsilon_{j} \ge \mathbb{E} \biggl( \int_{0}^{\tau} e^{-\lambda r_{j}^{m} s}
	\left(L(\widetilde{b}_{s}^{j}) + F_{j}(Y_{s}^{j}) \right) \, \mathrm{d}s + 
	e^{-\lambda r_{j}^{m} \tau} u_{j} (Y_{\tau}^{j}) \biggr).
\end{equation}
Indeed, subtracting the identity for the constant $u(0)$ and dividing by
$r_{j}^{\alpha}$ gives \eqref{e:scaled-stopped}; the control term is
invariant because $m = \alpha + \beta$ and $q(m - 1) = \beta$.

For $T > 0$ and $M > 1$, set $R := M T^{1 / m}$ and 
\begin{equation} \label{def:tau_j}
E_{j}(T) := \mathbb{E} \int_{0}^{T} \lvert \widetilde{b}_{s}^{j} \rvert^{q} \, \mathrm{d}s, \qquad \tau_{j} := T \wedge \inf\{s \ge 0 : \lvert Y_{s}^{j} - \xi \rvert \ge R\}.
\end{equation}

The next lemma bounds the energy of the rescaled controls and the probability of large excursions.
\begin{lem} \label{l:Et}
There exist $C_{E}, C_{X} > 0$, independent of $\xi$, $T$ and $M$, such that, for every $T > 0$,
\begin{equation} \label{e:Emain}
\limsup_{j \to \infty} E_{j}(T) \le C_{E} T^{\alpha / m},
\end{equation}
and, for every $M > 1$,
\begin{equation} \label{e:tmain}
\limsup_{j \to \infty} \mathbb{P}(\tau_{j} < T) \le C_{X} M^{-q}.
\end{equation}
\end{lem}
\begin{proof}
Write $L(b) = C_{L} \lvert b \rvert^{q}$, where $C_{L} := (p - 1) p^{-q}$. Since $f \ge -c_{f}$ and $J_{\lambda}(x_{j}; b^{j}) < \infty$, we have $E_{j}(T) < \infty$ for every $T > 0$. Moreover, $F_{j}(y) \ge - C r_{j}^{\beta}$ for $y \in \Omega_{j}$, where $C := c_{f} + \lvert \lambda u(0) \rvert$.
Taking $\tau = T$ in \eqref{e:scaled-stopped} and using \eqref{e:scaled-holder}, we obtain
\begin{align*}
	C_{L} e^{-\lambda r_{j}^{m} T} E_{j}(T) 
	&\le 
	u_{j}(\xi) + \varepsilon_{j} - 
	e^{-\lambda r_{j}^{m} T} \, \mathbb{E} \, u_{j}(Y_{T}^{j}) + 
	C r_{j}^{\beta} T
	\\ &\le 
	e^{-\lambda r_{j}^{m} T} \, \mathbb{E} \left( u_{j}(\xi) - u_{j}(Y_{T}^{j}) 
	\right) + \left(1 - e^{-\lambda r_{j}^{m} T} \right) u_{j}(\xi) + \varepsilon_{j} + o_{j}(1) 
	\\ &\le
	C_{H} \, \mathbb{E} \, \lvert Y_{T}^{j} - \xi \rvert^{\alpha} + o_{j}(1).
\end{align*} 
Then, for sufficiently large $j$,
\begin{equation} \label{e:Ea}
E_{j}(T) \le C_{HL} \, \mathbb{E} \, \lvert Y_{T}^{j} - \xi \rvert^{\alpha} + o_{j}(1),
\end{equation} 
where $C_{HL} := \frac{2 C_{H}}{C_{L}}$.
Since $Y_{T}^{j} - \xi = \int_{0}^{T} \widetilde{b}_{s}^{j} \, \mathrm{d}s + \sqrt{2 \varsigma_{j}} \, \widetilde{W}_{T}^{j}$ and $\alpha < 1$, by H\"older's inequality, we obtain
\begin{equation} \label{e:Ea2}
\mathbb{E} \, \lvert Y_{T}^{j} - \xi \rvert^{\alpha} \le T^{\alpha / p} E_{j}(T)^{\alpha / q} + (2 \varsigma_{j})^{\alpha / 2} \mathbb{E} \lvert \widetilde{W}_{T}^{j} \rvert^{\alpha} = T^{\alpha / p} E_{j}(T)^{\alpha / q} + C_{W} \varsigma_{j}^{\alpha / 2} T^{\alpha / 2},
\end{equation}
where $C_{W} := 2^{\alpha / 2} \mathbb{E} \lvert W_{1} \rvert^{\alpha}$.
By Young's inequality, we have
\begin{equation} \label{e:Ea3}
C_{HL} T^{\alpha / p} E_{j}(T)^{\alpha / q} \le \frac{\alpha}{q} E_{j}(T) + \Bigl( 1 - \frac{\alpha}{q} \Bigr) (C_{HL} T^{\alpha / p})^{\frac{q}{q - \alpha}}.
\end{equation} 
Moreover, $\frac{\alpha / p}{1 - \alpha / q} = \frac{\alpha}{m}$ and $\limsup_{j} \varsigma_{j}^{\alpha / 2} T^{\alpha / 2} \le T^{\alpha / m}$, with equality in the critical case and with limit $0$ in the supercritical case. Combining \eqref{e:Ea}, \eqref{e:Ea2} and \eqref{e:Ea3}, we obtain \eqref{e:Emain}.

By H\"older's inequality and \eqref{e:Emain}, using $q - 1 + \frac{\alpha}{m} = \frac{q}{m}$,
\[
\limsup_{j \to \infty} \, 
\mathbb{E} \left(\int_{0}^{T} \lvert \widetilde{b}_{s}^{j} \rvert \, \mathrm{d}s \right)^{q} 
\le C_{E} T^{q / m} = C_{E} \Bigl(\frac{R}{M} \Bigr)^{q}.
\]
On $\{\tau_{j} < T\}$ we have $R = \lvert Y_{\tau_{j}}^{j} - \xi \rvert \le \int_{0}^{T} \lvert \widetilde{b}_{s}^{j} \rvert \, \mathrm{d}s + \sqrt{2 \varsigma_{j}} \sup_{s \le T} \lvert \widetilde{W}_{s}^{j} \rvert$. Hence, by Chebyshev's and Doob's inequalities,
\[
\mathbb{P}(\tau_{j} < T) \le \frac{2^{q}}{R^{q}} \mathbb{E} \Bigl(\int_{0}^{T} \lvert \widetilde{b}_{s}^{j} \rvert \, \mathrm{d}s\Bigr)^{q} + \frac{32 N \varsigma_{j} T}{R^{2}}.
\]
The first term has $\limsup$ at most $2^{q} C_{E} M^{-q}$. The second equals $32 N \varsigma_{j} M^{-2} T^{(m - 2) / m}$; it is $32 N M^{-2} \le 32 N M^{-q}$ in the critical case, since $M > 1$ and $q < 2$, and it tends to $0$ in the supercritical case. This gives \eqref{e:tmain} with $C_{X} := 2^{q} C_{E} + 32 N$.
\end{proof}

\subsection{A scaled occupation estimate}
Under Condition~\ref{cond:supercritical}, the continuity of $a$ on $\overline{\Omega}_{\delta_{0}}$ and its positivity on $\partial\Omega$ allow us, after decreasing $\delta_{0}$ if necessary, to assume that $a \ge a_{*}$ on $\overline{\Omega}_{\delta_{0}}$ for some $a_{*} > 0$.

The next lemma establishes an occupation bound for near-optimal trajectories, needed in the comparison argument.

\begin{lem} \label{l:scaled-occupation}
For every $T > 0$ and $M > 1$,
\begin{equation} \label{e:scaled-occupation}
\limsup_{j \to \infty} \mathbb{E} \int_{0}^{\tau_{j}} d_{j}(Y_{s}^{j})^{-\beta} \, \mathrm{d}s < \infty.
\end{equation}
\end{lem}

\begin{proof}
For $s \le \tau_{j}$ we have $Y_{s}^{j} \in \Omega_{j} \cap \overline{B}(\xi, R)$, and for $y$ in this set $\operatorname{dist}(r_{j} y, \partial \Omega) \le r_{j} \lvert y \rvert \le r_{j} (\lvert \xi \rvert + R)$, since $0 \in \partial \Omega$. Hence, for all large $j$, $r_{j} y \in \Omega_{\delta_{0}}$, so $d(r_{j} y) = \operatorname{dist}(r_{j} y, \partial \Omega)$ and, for some $C_{d} > 0$,
\[
\lvert \nabla d_{j} \rvert = 1, \qquad \lvert \Delta d_{j} \rvert \le C_{d} r_{j}, \qquad 0 < d_{j} \le \lvert \xi \rvert + R \qquad \text{on } \Omega_{j} \cap \overline{B}(\xi, R).
\]

{\it Critical case}. Here $\beta = q$, $\alpha = 2 - q$, and $\varsigma_{j} = 1$. Using Young's inequality and the identity $p(1 - q) = -q$, we obtain, on $\Omega_{j} \cap \overline{B}(\xi, R)$, for every $b \in \mathbb{R}^{N}$,
\[
\begin{aligned}
(\Delta + b \cdot \nabla) d_{j}^{\alpha}
&= \alpha d_{j}^{1 - q} \bigl(\Delta d_{j} + b \cdot \nabla d_{j} \bigr) - \alpha (q - 1) d_{j}^{-q} 
\\ &\le 
-c_{0} d_{j}^{-q} + C_{0} \bigl(\lvert b \rvert^{q} + r_{j}^{q} \bigr),
\end{aligned}
\]
where $c_{0} := \alpha (q - 1) / p > 0$ and $C_{0} := 2^{q - 1} \alpha \max \{1, C_{d}^{q}\} / q$.

For $i \in \mathbb{N}$ set $\tau_{i,j} := \tau_{j} \wedge \inf \{ s \ge 0 : d_{j}(Y_{s}^{j}) \le i^{-1} \}$. Apply It\^{o}'s formula to $d_{j}^{\alpha}(Y_{s}^{j})$ on $[0, \tau_{i,j}]$, where $\nabla d_{j}^{\alpha}$ is bounded, take expectations, and use the last inequality with $b = \widetilde{b}_{s}^{j}$. Since $0 \le d_{j}^{\alpha} \le (\lvert \xi \rvert + R)^{\alpha}$ there,
\[
c_{0} \, \mathbb{E} \int_{0}^{\tau_{i, j}} d_{j}(Y_{s}^{j})^{-q} \, \mathrm{d}s \le (\lvert \xi \rvert + R)^{\alpha} + C_{0} E_{j}(T) + C_{0} r_{j}^{q} T.
\]
Since $d_{j}(Y^{j})$ is continuous and positive on the compact interval $[0, \tau_{j}]$, it is bounded below pathwise, so, almost surely, $\tau_{i, j} = \tau_{j}$ for all large $i$. Monotone convergence and \eqref{e:Emain} give \eqref{e:scaled-occupation}.

{\it Supercritical case}. By \eqref{e:scaled-forcing}, for all large $j$, we have $F_{j}(y) \ge a_{*} d_{j}(y)^{-\beta} - \lambda r_{j}^{\beta} \lvert u(0) \rvert$ on $\Omega_{j} \cap \overline{B}(\xi, R)$. Take $\tau = \tau_{j}$ in \eqref{e:scaled-stopped}. For large $j$, $e^{-\lambda r_{j}^{m} s} \ge 1 / 2$ for $s \le T$. Therefore
\[
\frac{a_{*}}{2} \, \mathbb{E} \int_{0}^{\tau_{j}} d_{j}(Y_{s}^{j})^{-\beta} \, \mathrm{d}s \le u_{j}(\xi) - \mathbb{E} \bigl(e^{-\lambda r_{j}^{m} \tau_{j}} u_{j}(Y_{\tau_{j}}^{j}) \bigr) + \varepsilon_{j} + \lambda r_{j}^{\beta} \lvert u(0) \rvert T.
\]
Since $\lvert Y_{\tau_{j}}^{j} - \xi \rvert \le R$ and $\tau_{j} \le T$, \eqref{e:scaled-holder} and \eqref{e:scaled-growth} give
\[
u_{j}(\xi) - \mathbb{E} \bigl(e^{-\lambda r_{j}^{m} \tau_{j}} u_{j}(Y_{\tau_{j}}^{j}) \bigr) \le C_{H} R^{\alpha} + \bigl(1 - e^{-\lambda r_{j}^{m} T} \bigr) C_{H} \lvert \xi \rvert^{\alpha}.
\]
Letting $j \to \infty$ gives $\limsup_{j} \mathbb{E} \int_{0}^{\tau_{j}} d_{j}(Y_{s}^{j})^{-\beta} \, \mathrm{d}s \le 2 C_{H} R^{\alpha} / a_{*}$, which proves \eqref{e:scaled-occupation}.
\end{proof}

\subsection{Proof of tangential invariance}

We show that every blow-up limit $\Phi$ is invariant under tangential translations by constructing, from a near-optimal trajectory starting at $\xi$, a second trajectory starting near $\xi + h$. We introduce an additional tangential drift to gradually cancel the shift and a small inward correction to preserve the state constraint up to a stopping time. We then use the dynamic programming inequalities to compare their costs and show that the comparison error vanishes in the limit.
\begin{prop}\label{p:ti}
Let $\Phi$ be a blow-up limit as in \eqref{e:blowup-limit}. Then there is $\varphi \in C([0, \infty))$ such that
\begin{equation}\label{e:one-dimensional}
\Phi(z, \rho) = \varphi(\rho), \qquad (z, \rho) \in \overline{H}, \qquad \varphi(0) = 0.
\end{equation}
\end{prop}
\begin{proof}
{\it Step 1: tangential shift and comparison inequality}.
Fix $T > 0$, $M > 1$, $\xi \in H$, $h = (h^{\prime}, 0) \in \mathbb{R}^{N - 1} \times \{0\}$, and let $\varepsilon_{j}$, $b^{j}$, $Y^{j}$, $R$ and $\tau_{j}$ be as in Section \ref{ss:snoc} for this $\xi$.
Let $h_{s} = (h^{\prime}_{s}, 0) := \left(1 - s / T \right) h$ for $s \in [0, T]$. 
By \eqref{def:tau_j}, $Y_{s}^{j} \in \Omega_{j} \cap \overline{B}(\xi, R)$ for $s \in [0, \tau_{j}]$. Set
\begin{equation} \label{def:K}
\mathcal{K} := \overline{B}(\xi, R) + \left\{ \zeta_{1} h + \zeta_{2} e_{N} \mid \zeta_{1}, \zeta_{2} \in [0, 1] \right\}.
\end{equation} 
Choose a bounded cylinder $\mathcal{C}$ centered at $0$ with $\overline{B}(0, 2 \sup_{y \in \mathcal{K}} \lvert y \rvert) \subset \mathcal{C}$. Since $0 \in \partial \Omega_{j}$, every nearest boundary point of a point $y \in \mathcal{K} \cap \Omega_{j}$ lies in $\overline{B}(0, 2 \lvert y \rvert) \subset \mathcal{C}$, and, by Lemma~\ref{l:blowup}, $\Omega_{j} \cap \mathcal{C} = \{(z, \rho) \in \mathcal{C} : \rho > g_{j}(z)\}$ for large $j$.
Let $\mathcal{C}_{\mathrm{tan}}$ be the base of $\mathcal{C}$, and let $\ell_{j} := \sup_{z \in \mathcal{C}_{\mathrm{tan}}} \lvert \nabla g_{j}(z) \rvert$. By \eqref{e:scaled-graph}, $\ell_{j} = O(r_{j})$. Set $\delta_{j} := \ell_{j} \lvert h \rvert$, so that $\delta_{j} \to 0$ as $j \to \infty$. Define
\begin{equation*}
\widehat{Y}_{s}^{j} := Y_{s}^{j} + h_{s} + \delta_{j} e_{N}, \qquad 
\widehat{b}_{s}^{j} := \widetilde{b}_{s}^{j} - \frac{h}{T}.
\end{equation*}
Then
\[
\mathrm{d} \widehat{Y}_{s}^{j} = \widehat{b}_{s}^{j} \, \mathrm{d}s + \sqrt{2 \varsigma_{j}} \, \mathrm{d} \widetilde{W}_{s}^{j},
 \qquad \widehat{Y}_{0}^{j} = \xi + h + \delta_{j} e_{N}.
\] 
Write $Y_{s}^{j} = (Z_{s}^{j}, P_{s}^{j}) \in \mathbb{R}^{N - 1} \times \mathbb{R}$.
By \eqref{def:K}, since $\delta_{j} \le 1$ for large $j$, we have $\widehat{Y}_{s}^{j} \in \mathcal{K}$ for $s \in [0, \tau_{j}]$. Moreover, since $g_{j}(Z_{s}^{j} + h_{s}^{\prime}) - g_{j}(Z_{s}^{j}) \le \ell_{j} \lvert h_{s}^{\prime} \rvert \le \delta_{j}$, the graph representation gives $\widehat{Y}_{s}^{j} \in \Omega_{j}$ for $s \in [0, \tau_{j}]$. The processes $Y^{j}$ and $\widehat{Y}^{j}$ have the same natural filtration because their difference is deterministic, so $\tau_{j}$ is also a stopping time for $\widehat{Y}^{j}$.

Apply the verification inequality \eqref{e:uu} to the rescaled process $\widehat{X}_{t}^{j} = r_{j} \widehat{Y}_{t / r_{j}^{m}}^{j}$ up to the stopping time $r_{j}^{m} \tau_{j}$. Scaling this inequality and comparing it with \eqref{e:scaled-stopped} at $\tau_{j}$, we obtain
\begin{align} \nonumber
u_{j}(\xi + h + \delta_{j} e_{N}) - u_{j}(\xi) - \varepsilon_{j}
 \le {}&
\mathbb{E} \int_{0}^{\tau_{j}} e^{-\lambda r_{j}^{m} s} 
\left| L\left(\widetilde{b}_{s}^{j} - \frac{h}{T} \right) - L(\widetilde{b}_{s}^{j}) 
\right| \, \mathrm{d}s
\\ \nonumber &+ 
\mathbb{E} \int_{0}^{\tau_{j}} e^{-\lambda r_{j}^{m} s}
\bigl(F_{j}(\widehat{Y}_{s}^{j}) - F_{j}(Y_{s}^{j}) \bigr) \, \mathrm{d}s
\\ \label{e:prelimit-comparison} &+ 
\mathbb{E} \left\lvert u_{j}(\widehat{Y}_{\tau_{j}}^{j}) - u_{j}(Y_{\tau_{j}}^{j}) \right\rvert.
\end{align}

Since $L(b) = C_{L} \lvert b \rvert^{q}$ with $q \in (1, 2)$, there is $C > 0$ such that $\lvert L(v + w) - L(v) \rvert \le C \bigl(\lvert v \rvert^{q - 1} \lvert w \rvert + \lvert w \rvert^{q} \bigr)$ for all $v, w \in \mathbb{R}^{N}$. Moreover, by H\"{o}lder's inequality, $\mathbb{E} \int_{0}^{T} \lvert \widetilde{b}_{s}^{j} \rvert^{q - 1} \, \mathrm{d}s \le T^{1 / q} E_{j}(T)^{(q - 1) / q}$. Combining these inequalities with \eqref{e:Emain} and using $-1 + 1 / q + \alpha(q - 1) / (mq) = -\theta / m$, we obtain, after enlarging $C$,
\begin{equation} \label{e:control-error}
\limsup_{j \to \infty} \mathbb{E} \int_{0}^{\tau_{j}} e^{-\lambda r_{j}^{m} s}
\left|
L\left( \widetilde{b}_{s}^{j} - \frac{h}{T} \right) - L(\widetilde{b}_{s}^{j})
\right| \, \mathrm{d}s
 \le C \left( \lvert h \rvert T^{-\theta / m} + \lvert h \rvert^{q} T^{1 - q} \right).
\end{equation}

Since $h_{T} = 0$, on $\{\tau_{j} = T\}$ we have $\widehat{Y}_{T}^{j} - Y_{T}^{j} = \delta_{j} e_{N}$, whereas on $\{\tau_{j} < T\}$ we have $\lvert \widehat{Y}_{\tau_{j}}^{j} - Y_{\tau_{j}}^{j} \rvert 
 \le \lvert h \rvert + \delta_{j}$. Hence \eqref{e:scaled-holder} gives
\[
\mathbb{E} \left|
u_{j}(\widehat{Y}_{\tau_{j}}^{j}) - u_{j}(Y_{\tau_{j}}^{j})
\right|
 \le C_{H} \delta_{j}^{\alpha} + C_{H}(\lvert h \rvert + \delta_{j})^{\alpha} \, \mathbb{P}(\tau_{j} < T).
\]
By \eqref{e:tmain}, letting $j \to \infty$ gives
\begin{equation} \label{e:terminal-error}
\limsup_{j \to \infty} \mathbb{E} \, \left\lvert u_{j}(\widehat{Y}_{\tau_{j}}^{j}) - 
u_{j}(Y_{\tau_{j}}^{j}) \right\rvert \le C_{X} C_{H} \lvert h \rvert^{\alpha} M^{-q}.
\end{equation}
Also, by \eqref{e:scaled-holder},
\begin{equation} \label{e:uci}
\lvert u_{j}(\xi + h + \delta_{j} e_{N}) - u_{j}(\xi + h) \rvert \le C_{H} \delta_{j}^{\alpha} \to 0 \quad \text{as } j \to \infty.
\end{equation}

It remains to show that the state cost term in \eqref{e:prelimit-comparison} has nonpositive $\limsup$.

{\it Step 2: state cost error}.
For all sufficiently large $j$ and all $s \in [0, \tau_{j}]$, both $r_{j} Y_{s}^{j}$ and $r_{j} \widehat{Y}_{s}^{j}$ lie in $\Omega_{\delta_{0}}$. Put
\[
\begin{aligned}
d_{s}^{j} &:= d_{j}(Y_{s}^{j}), &\qquad \widehat{d}_{s}^{j} &:= d_{j}(\widehat{Y}_{s}^{j}), \\
a_{s}^{j} &:= a(r_{j} Y_{s}^{j}), &\qquad \widehat{a}_{s}^{j} &:= a(r_{j} \widehat{Y}_{s}^{j}).
\end{aligned}
\]
Since the rescaled state cost involves the singular factor $d_{j}^{-\beta}$, we compare the distances to the boundary along the two paths. For large $j$, $d_{j} = \operatorname{dist}(\cdot, \partial \Omega_{j})$ on $\Omega_{j} \cap \mathcal{C}$. Let $y = (z, \rho) \in \mathcal{K} \cap \Omega_{j}$. The vertical segment gives $d_{j}(y) \le \rho - g_{j}(z)$. If $(w, g_{j}(w))$ is a nearest boundary point of $y$, it lies in $\mathcal{C}$ by Step~1, and since the base $\mathcal{C}_{\mathrm{tan}}$ is convex,
\[
d_{j}(y) \le \rho - g_{j}(z) \le \lvert \rho - g_{j}(w) \rvert + \ell_{j} \lvert w - z \rvert \le \sqrt{1 + \ell_{j}^{2}} \, d_{j}(y).
\]
Applying the rightmost inequality to $\widehat{Y}_{s}^{j} = (Z_{s}^{j} + h_{s}^{\prime}, P_{s}^{j} + \delta_{j})$ and the leftmost inequality to $Y_{s}^{j}$, and using $P_{s}^{j} + \delta_{j} - g_{j}(Z_{s}^{j} + h_{s}^{\prime}) \ge P_{s}^{j} - g_{j}(Z_{s}^{j})$ by the definition of $\delta_{j}$, we obtain $\widehat{d}_{s}^{j} \ge d_{s}^{j} / \sqrt{1 + \ell_{j}^{2}}$ for $s \in [0, \tau_{j}]$, that is,
\begin{equation*}
(\widehat{d}_{s}^{j})^{-\beta} \le (1 + \gamma_{j}) \, (d_{s}^{j})^{-\beta}, \qquad s \in [0, \tau_{j}], \qquad \gamma_{j} := (1 + \ell_{j}^{2})^{\beta / 2} - 1.
\end{equation*}
Since $\ell_{j} = O(r_{j})$, $\gamma_{j} \to 0$ as $j \to \infty$.

Put $a_{0} := a(0)$ and set $\omega_{j} := \operatorname*{ess\,sup}_{y \in \mathcal{K} \cap \Omega_{j}} \lvert a(r_{j}y) - a_{0} \rvert$. By the assumptions on $a$ and its boundary trace, $\omega_{j} \to 0$ as $j \to \infty$. The constant $-\lambda r_{j}^{\beta}u(0)$ cancels, and since $a_{0} \ge 0$,
\[
\begin{aligned}
F_{j}(\widehat{Y}_{s}^{j}) - F_{j}(Y_{s}^{j})
&=
a_{0} \bigl((\widehat{d}_{s}^{j})^{-\beta} - (d_{s}^{j})^{-\beta}\bigr) + (\widehat{a}_{s}^{j} - a_{0})(\widehat{d}_{s}^{j})^{-\beta} - (a_{s}^{j}-a_{0})(d_{s}^{j})^{-\beta}
\\ &\le
\bigl(a_{0} \gamma_{j} + (2 + \gamma_{j}) \, \omega_{j} \bigr) (d_{s}^{j})^{-\beta}.
\end{aligned}
\] 
For measurable $a$, this inequality holds for a.e. $s \in [0,\tau_{j}]$ almost surely by the occupation estimate of \cite[Th.~4.1]{WLW23} for the two controlled diffusions.
By Lemma~\ref{l:scaled-occupation}, since $\gamma_{j},\omega_{j}\to0$,
\begin{equation} \label{e:forcing-error}
\limsup_{j \to \infty} \mathbb{E} \int_{0}^{\tau_{j}} e^{-\lambda r_{j}^{m}s} \bigl(F_{j}(\widehat{Y}_{s}^{j}) - F_{j}(Y_{s}^{j}) \bigr) \, \mathrm{d}s \le 0.
\end{equation}

{\it Step 3: conclusion}.
Combining \eqref{e:prelimit-comparison} with \eqref{e:control-error}, \eqref{e:terminal-error}, \eqref{e:uci}, and \eqref{e:forcing-error} and letting first $j \to \infty$, then $M \to \infty$, and finally $T \to \infty$, gives $\Phi(\xi + h) \le \Phi(\xi)$. Repeating the argument from $\xi + h$ with tangential shift $-h$ gives the reverse inequality. Therefore
\begin{equation} \label{e:tangential-invariance}
\Phi(\xi + h) = \Phi(\xi).
\end{equation}
Since $\xi \in H$ and $h$ are arbitrary and $\Phi$ is continuous on $\overline{H}$, \eqref{e:tangential-invariance} holds for all $\xi \in \overline{H}$. Hence $\varphi(\rho) := \Phi(0, \rho)$ is continuous on $[0, \infty)$, satisfies $\varphi(0) = \Phi(0) = 0$, and gives \eqref{e:one-dimensional}.
\end{proof}

\section{Classification of blow-up limits} \label{s:lo}
By Lemma~\ref{l:blowup}, after passing to a subsequence, $u_{j} \to \Phi$ locally uniformly in $H$ for some $\Phi \in C(\overline{H})$. In this section we prove that $\Phi$ is one-dimensional and uniquely determined, identify it explicitly, and show that the convergence holds in $C^{1}_{\mathrm{loc}}(H)$.

\subsection{The rescaled equation}
By \eqref{e:u} and \eqref{e:LL-regularity}, the rescaled functions satisfy
\begin{equation} \label{e:scaled-PDE}
- \varsigma_{j} \Delta u_{j} + \lvert \nabla u_{j} \rvert^{p} + \lambda r_{j}^{m} u_{j} = F_{j} \qquad \text{a.e. in } \Omega_{j}.
\end{equation}
By \eqref{e:LL-gradient},
\begin{equation} \label{e:scaled-gradient}
\lvert \nabla u_{j}(y) \rvert \le C_{G} \, d_{j}(y)^{-\theta}, \qquad y \in \Omega_{j}.
\end{equation}

Recall that $a_{0} := a(0)$. Under Condition~\ref{cond:critical} or \ref{cond:supercritical}, the assumptions on $a$ give $a(r_{j} y) \to a_{0}$ in $L^{\infty}_{\mathrm{loc}}(H)$, while Lemma~\ref{l:blowup} gives $d_{j}(y) \to \rho$ uniformly on every compact subset of $H$. On each such subset, $r_{j} y \in \Omega_{\delta_{0}}$ for large $j$, and hence
\begin{equation} \label{e:forcing-limit}
F_{j}(y) = a(r_{j} y) \, d_{j}(y)^{-\beta} - \lambda r_{j}^{\beta} u(0) \longrightarrow a_{0} \rho^{-\beta} \qquad \text{in } L_{\mathrm{loc}}^{\infty}(H).
\end{equation}

\subsection{Boundary growth estimate}
We establish a uniform version of the boundary growth property of $u$ obtained by Lasry and Lions \cite[Th.~IV.2, Cor.~IV.3]{LL89}. This estimate passes to blow-up limits and will be used to select the limiting profiles.
\begin{lem} \label{l:asn}
There exist $\varrho_{0}, \eta_{0} > 0$ such that
\begin{equation} \label{e:asn}
\inf_{\overline{\Omega} \cap \overline{B}(x_{0}, \varrho)} u \le u(x_{0}) - \eta_{0} \varrho^{\alpha}, \qquad x_{0} \in \partial \Omega, \quad \varrho \in (0, \varrho_{0}].
\end{equation}
Consequently, every blow-up limit $\Phi$ as in \eqref{e:blowup-limit} satisfies
\begin{equation} \label{e:pa}
\inf_{\overline{H} \cap \overline{B}(0, \varrho)} \Phi \le -\eta_{0} \varrho^{\alpha}, \qquad \varrho > 0.
\end{equation}
\end{lem}

\begin{proof}
{\it Step 1: prelimit estimate}.
Recall that the maximal solution $u$ satisfies the state-constraint boundary condition \eqref{e:state-constraint}.

Fix $x_{0} \in \partial \Omega$ and $\varrho \in (0, \varrho_{0}]$, where $\varrho_{0} \in (0, \delta_{0})$ will be chosen below. For $\varepsilon > 0$, define
\[
\psi(x) := u(x_{0}) + \varepsilon \varrho^{\alpha - 2} \left( d(x)(d(x) - \varrho) - \frac{\lvert x - x_{0} \rvert^{2}}{2N}\right), \qquad x \in \overline{\Omega}.
\]
Set $U := \Omega \cap B(x_{0}, \varrho)$. Using $0 \le d \le \varrho$, $\lvert \nabla d \rvert = 1$, and the boundedness of $\Delta d$ in $\Omega_{\delta_{0}}$, we obtain $\lvert \nabla \psi \rvert \le 2 \varepsilon \varrho^{\alpha - 1}$ and $\Delta \psi \ge \frac{\varepsilon}{2} \varrho^{\alpha - 2}$ in $U$, for sufficiently small $\varrho_{0}$.

Under Condition~\ref{cond:critical}, we have $\alpha - 2 = -\beta$ and $f \ge -c_{f}$, while under Condition~\ref{cond:supercritical}, we have $f \ge a_{*} \varrho^{-\beta}$ and $\Delta \psi \ge 0$. Thus, in both cases, for sufficiently small $\varepsilon > 0$, we have $f + \Delta \psi \ge \frac{\varepsilon}{2} \varrho^{-\beta} - c_{f}$ a.e. in $U$. Since $p(\alpha - 1) = -\beta$, we deduce that
\[
-\Delta \psi + \lvert \nabla \psi \rvert^{p} + \lambda u - f \le \left(-\frac{\varepsilon}{2} + 2^{p} \varepsilon^{p}\right) \varrho^{-\beta} + \lambda \lVert u \rVert_{\infty} + c_{f}.
\]
Choosing $\varepsilon$ small and then decreasing $\varrho_{0}$, independently of $x_{0}$, we obtain
\begin{equation} \label{e:growth-strict}
-\Delta \psi + \lvert \nabla \psi \rvert^{p} + \lambda u \le f - 1 \qquad \text{a.e. in } U.
\end{equation}

Let $x_{*} \in \overline{U}$ be a minimum point of $u - \psi$. By the boundary condition \eqref{e:state-constraint}, $x_{*} \notin \partial \Omega \cap B(x_{0}, \varrho)$. If $x_{*} \in U$, then $\nabla u(x_{*}) = \nabla \psi(x_{*})$. Combining \eqref{e:u} with \eqref{e:growth-strict} and using the continuity of the gradients, we obtain $-\Delta(u - \psi) \ge 1 / 2$ a.e. in a neighborhood of $x_{*}$, contradicting the strong minimum principle \cite[Th.~9.6]{GT01}. Consequently, $x_{*} \in \overline{U} \cap \partial B(x_{0}, \varrho)$. Since $u(x_{0}) = \psi(x_{0})$, we have $u(x_{*}) \le \psi(x_{*})$. Using $\lvert x_{*} - x_{0} \rvert = \varrho$ and $0 \le d(x_{*}) \le \varrho$, we obtain
\[
u(x_{*}) \le \psi(x_{*}) \le u(x_{0}) - \frac{\varepsilon}{2N} \varrho^{\alpha}.
\]
This proves \eqref{e:asn} with $\eta_{0} := \varepsilon / (2N)$.

{\it Step 2: blow-up limits}. Finally, applying \eqref{e:asn} at the physical scale $r_{j} \varrho$ and dividing by $r_{j}^{\alpha}$ gives, for every fixed $\varrho > 0$ and all large $j$,
\[
	\inf_{\overline{\Omega}_{j} \cap \overline{B}(0, \varrho)} u_{j} \le -\eta_{0} \varrho^{\alpha}.
\]
Let $\delta > 0$ and choose $y_{j} \in \overline{\Omega}_{j} \cap \overline{B}(0, \varrho)$ with $u_{j}(y_{j}) \le - \eta_{0} \varrho^{\alpha}$. By \eqref{e:scaled-graph}, for all large $j$, $y_{j} + \delta e_{N}$ lies in the compact set $K_{\delta} := \overline{B}(0, \varrho + \delta) \cap \{\rho \ge \delta / 2\} \subset H \cap \Omega_{j}$, and \eqref{e:scaled-holder} gives $u_{j}(y_{j} + \delta e_{N}) \le - \eta_{0} \varrho^{\alpha} + C_{H} \delta^{\alpha}$. Since $u_{j} \to \Phi$ uniformly on $K_{\delta}$, we obtain $\min_{\overline{H} \cap \overline{B}(0, \varrho + \delta)} \Phi \le - \eta_{0} \varrho^{\alpha} + C_{H} \delta^{\alpha}$, and letting $\delta \downarrow 0$ gives \eqref{e:pa}.
\end{proof}

\subsection{The critical case}
\begin{lem} \label{l:critical-blowup}
Assume Condition~\ref{cond:critical}. Then $u_{j} \to \Phi$ in $C^{1}_{\mathrm{loc}}(H)$, and
\begin{equation} \label{e:critical-unique-blowup}
\Phi(z, \rho) = -\frac{A(a_{0})}{\alpha} \rho^{\alpha}, \qquad \nabla \Phi(z, \rho) = -A(a_{0}) \rho^{-\theta} e_{N}, \qquad (z, \rho) \in H.
\end{equation}
\end{lem}

\begin{proof}
{\it Step 1: $C^{1}$ convergence and the limit equation}. Here $\varsigma_{j} = 1$. Fix $k > N$ and a compact set $K \Subset H$ that is the closure of a smooth open set $K^{\circ}$, and choose $\chi \in C_{c}^{\infty}(H)$ with $0 \le \chi \le 1$ and $\chi \equiv 1$ near $K$. By Lemma~\ref{l:blowup}, there exists $j_{K} \in \mathbb{N}$ such that $\operatorname{supp} \chi \subset \Omega_{j}$ for all $j \ge j_{K}$. By the Calder\'{o}n--Zygmund estimate \cite[Cor.~9.10]{GT01}, there exists $C_{N, k} > 0$ such that
\[
\lVert D^{2} u_{j} \rVert_{L^{k}(K)} \le \lVert D^{2} (\chi u_{j}) \rVert_{L^{k}(H)} \le C_{N, k} \lVert \Delta (\chi u_{j}) \rVert_{L^{k}(H)}.
\]
Moreover, for some $C_{\chi} > 0$, $\lVert \Delta(\chi u_{j}) \rVert_{L^{k}(H)} \le C_{\chi} \bigl\lVert \lvert \Delta u_{j} \rvert + \lvert \nabla u_{j} \rvert + \lvert u_{j} \rvert \bigr\rVert_{L^{k}(\operatorname{supp} \chi)}$, and $\Delta u_{j}$, $\nabla u_{j}$ and $u_{j}$ are essentially bounded on $\operatorname{supp} \chi$, uniformly in $j$, by \eqref{e:scaled-PDE}, \eqref{e:forcing-limit}, \eqref{e:scaled-gradient} and \eqref{e:scaled-growth}. Hence
\[
\sup_{j \ge j_{K}} \lVert u_{j} \rVert_{W^{2, k}(K^{\circ})} < \infty.
\]
Since $W^{2, k}(K^{\circ})$ embeds compactly into $C^{1}(K)$ for $k > N$, every subsequence of $(u_{j})$ has a further subsequence converging in $C^{1}(K)$, with $D^{2} u_{j}$ converging weakly in $L^{k}(K)$; by \eqref{e:blowup-limit}, the limit is $\Phi$. Since $K$ is arbitrary, $u_{j} \to \Phi$ in $C^{1}_{\mathrm{loc}}(H)$, $\Phi \in W^{2, k}_{\mathrm{loc}}(H)$, and passing to the limit in \eqref{e:scaled-PDE} gives
\begin{equation*}
-\Delta \Phi + \lvert \nabla \Phi \rvert^{p} = a_{0} \rho^{-q} \qquad \text{a.e. in } H.
\end{equation*}

{\it Step 2: one-dimensional reduction and classification}. By Proposition~\ref{p:ti}, $\nabla \Phi = \varphi^{\prime}(\rho) e_{N}$ and $\Delta \Phi = \varphi^{\prime \prime}(\rho)$, so
\begin{equation} \label{e:critical-limit-ODE}
-\varphi^{\prime \prime} + \lvert \varphi^{\prime} \rvert^{p} = a_{0} \rho^{-q}, \qquad \rho \in (0, \infty), \qquad \varphi(0) = 0.
\end{equation}
Since $\varphi \in W^{2, k}_{\mathrm{loc}}((0, \infty))$ and $\varphi^{\prime \prime} = \lvert \varphi^{\prime} \rvert^{p} - a_{0} \rho^{-q}$ is continuous, $\varphi \in C^{2}((0, \infty))$. Since $d_{j} \to \rho$ and $\nabla u_{j} \to \nabla \Phi$ locally uniformly, \eqref{e:scaled-gradient} passes to the limit:
\begin{equation} \label{e:php}
\lvert \varphi^{\prime}(\rho) \rvert \le C_{G} \rho^{-\theta}, \qquad \rho > 0.
\end{equation}
Set $w(\widetilde{\rho}) := -e^{\theta \widetilde{\rho}} \varphi^{\prime}(e^{\widetilde{\rho}})$ for $\widetilde{\rho} \in \mathbb{R}$. By \eqref{e:php}, $w$ is bounded, and, since $\theta + 1 = q$ and $p \theta = q$, \eqref{e:critical-limit-ODE} becomes
\begin{equation} \label{e:critical-autonomous}
w^{\prime} = a_{0} + \theta w - \lvert w \rvert^{p} =: S(w).
\end{equation}
Since $S$ is concave and $S(0) = a_{0} \ge 0$, it has exactly two zeros $B \le 0 < A(a_{0})$, where $A(a_{0})$ is the root in \eqref{e:critical-algebraic-B}, and $B = 0$ if and only if $a_{0} = 0$. Moreover, $S > 0$ on $(B, A(a_{0}))$ and $S < 0$ outside $[B, A(a_{0})]$. Consequently, every bounded solution of \eqref{e:critical-autonomous} on $\mathbb{R}$ other than $w \equiv A(a_{0})$ satisfies $w(\widetilde{\rho}) \to B$ as $\widetilde{\rho} \to -\infty$. Since $w$ is bounded, $\theta < 1$ and $\varphi(0) = 0$,
\begin{equation} \label{e:critical-profile-integral}
\varphi(\rho) = -\int_{0}^{\rho} w(\log s) \, s^{-\theta} \, \mathrm{d}s, \qquad \rho > 0.
\end{equation}
By Proposition~\ref{p:ti}, \eqref{e:pa} reads $\min_{[0, \rho]} \varphi \le -\eta_{0} \rho^{\alpha}$ for $\rho > 0$. Suppose that $w \not\equiv A(a_{0})$. If $B < 0$, then, by dominated convergence, $\varphi(\rho) / \rho^{\alpha} \to -B / \alpha > 0$ as $\rho \downarrow 0$, so $\varphi \ge 0$ near $0$, contradicting \eqref{e:pa}. If $B = 0$, then $w \ge 0$ and $w(\widetilde{\rho}) \to 0$ as $\widetilde{\rho} \to -\infty$, so $\varphi \le 0$ and $\varphi(\rho) = o(\rho^{\alpha})$ as $\rho \downarrow 0$; hence $\min_{[0, \rho]} \varphi = o(\rho^{\alpha})$ as $\rho \downarrow 0$, again contradicting \eqref{e:pa}. Therefore $w \equiv A(a_{0})$, and \eqref{e:critical-profile-integral} gives \eqref{e:critical-unique-blowup}.

Since every subsequence of $(u_{j})$ has a further subsequence converging to the same profile in $C^{1}_{\mathrm{loc}}(H)$, the full sequence converges in $C^{1}_{\mathrm{loc}}(H)$.
\end{proof}

\subsection{The supercritical case}
Here $\varsigma_{j} \to 0$, so the elliptic regularity argument from the critical case is no longer available. We instead identify the limiting profile using viscosity solutions.

{\it Viscosity solutions}. Let $\mathcal{O} \subset \mathbb{R}^{N}$ be open, let $\mathcal{S}(N)$ denote the space of symmetric $N \times N$ matrices, and let $\mathcal{G} \colon \mathcal{O} \times \mathbb{R} \times \mathbb{R}^{N} \times \mathcal{S}(N) \to \mathbb{R}$ be continuous. A function $v \in C(\mathcal{O})$ is a {\it viscosity subsolution} (resp.\ {\it supersolution}) of $\mathcal{G}(y, v, \nabla v, D^{2} v) = 0$ in $\mathcal{O}$ if, for every $\psi \in C^{2}(\mathcal{O})$ and every local maximum (resp.\ minimum) point $y_{0} \in \mathcal{O}$ of $v - \psi$,
\[
\mathcal{G}\bigl(y_{0}, v(y_{0}), \nabla \psi(y_{0}), D^{2} \psi(y_{0}) \bigr) \le 0 \qquad (\text{resp.\ } \ge 0),
\]
and a {\it viscosity solution} if it is both. For continuous functions, this definition is equivalent to the semijet formulation in \cite[Def.~2.2, Rem.~2.3]{CIL92}. For first-order equations, $C^{1}$ test functions suffice; see \cite[Ch.~II, Def.~1.1]{BCD97}.

\begin{lem} \label{l:super-profile}
Assume Condition~\ref{cond:supercritical}. Then
\begin{equation} \label{e:super-unique-blowup}
\Phi(z, \rho) = -\frac{a_{0}^{1 / p}}{\alpha} \rho^{\alpha}, \qquad \nabla \Phi(z, \rho) = -a_{0}^{1 / p} \rho^{-\theta} e_{N}, \qquad (z, \rho) \in H.
\end{equation}
\end{lem}

\begin{proof}
{\it Step 1: the limit equation}. By \eqref{e:scaled-gradient}, the functions $u_{j}$ are uniformly Lipschitz on compact subsets of $H$, and hence so is $\Phi$. Since $u_{j} \in W^{2, N}_{\mathrm{loc}}(\Omega_{j})$ by \eqref{e:LL-regularity} and $F_{j}$ is continuous, $u_{j}$ is a viscosity solution of \eqref{e:scaled-PDE} in $\Omega_{j}$ by \cite[Th.~I.2]{pL83}. As $j \to \infty$, we have $u_{j} \to \Phi$ and $F_{j} \to a_{0} \rho^{-\beta}$ locally uniformly in $H$, while $\varsigma_{j} \to 0$ and $\lambda r_{j}^{m} \to 0$. The standard stability argument for viscosity solutions, applied at local extrema of $u_{j} - \psi$ for smooth test functions $\psi$, yields
\begin{equation*}
\lvert \nabla \Phi \rvert^{p} = a_{0} \rho^{-\beta} \qquad \text{in } H
\end{equation*}
in the viscosity sense, cf. \cite[Sect.~6]{CIL92}. Hence, using Proposition~\ref{p:ti}, $\varphi$ is a viscosity solution of
\begin{equation} \label{e:super-limit-ODE}
\lvert \varphi^{\prime} \rvert^{p} = a_{0} \rho^{-\beta} \qquad \text{in } (0, \infty).
\end{equation}

{\it Step 2: classification}. Suppose that $\varphi$ has a local minimum at some $\rho_{0} > 0$. Taking the constant test function $\psi \equiv \varphi(\rho_{0})$, the supersolution property gives $0 \ge a_{0} \rho_{0}^{-\beta}$, a contradiction. Thus there is $\rho_{*} \in [0, \infty]$ such that $\varphi$ is nondecreasing on $(0, \rho_{*})$ and nonincreasing on $(\rho_{*}, \infty)$. By Proposition~\ref{p:ti}, \eqref{e:pa} reads $\min_{[0, \rho]} \varphi \le -\eta_{0} \rho^{\alpha}$ for $\rho > 0$. If $\rho_{*} > 0$, then, since $\varphi(0) = 0$ and $\varphi$ is continuous, $\varphi \ge 0$ on $[0, \rho_{*})$, contradicting \eqref{e:pa} for $\rho \in (0, \rho_{*})$. Hence $\rho_{*} = 0$, and $\varphi$ is nonincreasing on $(0, \infty)$.

By \cite[Ch.~II, Prop.~1.9]{BCD97}, $\varphi$ satisfies \eqref{e:super-limit-ODE} at every point of differentiability, hence almost everywhere. As $\varphi$ is nonincreasing and $\beta / p = \theta$, we obtain $\varphi^{\prime}(\rho) = -a_{0}^{1 / p} \rho^{-\theta}$ for a.e.\ $\rho > 0$. Integrating and using $\varphi(0) = 0$ and $\theta < 1$ yields $\varphi(\rho) = -a_{0}^{1 / p} \rho^{\alpha} / \alpha$, which proves \eqref{e:super-unique-blowup}.
\end{proof}

\subsection{Gradient convergence in the vanishing-viscosity limit}
To complete the analysis of the supercritical case, we prove a general gradient-convergence lemma for the vanishing-viscosity limit. We first establish the local regularity needed for its proof.
\begin{lem} \label{l:w3s}
Let $U \subset \mathbb{R}^{N}$ be open, $\varsigma > 0$, $\lambda \ge 0$ and $F \in W^{1, \infty}_{\mathrm{loc}}(U)$. Suppose that $v \in W^{2, k}_{\mathrm{loc}}(U)$ for all $k \in (1, \infty)$ and
\begin{equation} \label{e:w-equation}
-\varsigma \Delta v + \lvert \nabla v \rvert^{p} + \lambda v = F \qquad \text{a.e. in } U.
\end{equation}
Then $v \in W^{3, k}_{\mathrm{loc}}(U)$ for every $k \in (1, \infty)$, and, for $i = 1, \dots, N$,
\begin{equation} \label{e:w-differentiated}
-\varsigma \Delta(\partial_{i} v) + p \lvert \nabla v \rvert^{p - 2} \nabla v \cdot \nabla(\partial_{i} v) + \lambda \partial_{i} v = \partial_{i} F \qquad \text{a.e. in } U.
\end{equation}
\end{lem}
\begin{proof}
Fix $k \in (1, \infty)$. By Sobolev embedding and the assumed regularity of $v$, $\nabla v$ is continuous. Let $(\nabla v)_{\varepsilon}$ denote mollifications of $\nabla v$. On compact subsets of $U$, as $\varepsilon \to 0$, $(\nabla v)_{\varepsilon} \to \nabla v$ uniformly and in $W^{1, k}$, so, since the map $\lvert \cdot \rvert^{p}$ is $C^{1}$,
\[
\partial_{i} \bigl[\lvert (\nabla v)_{\varepsilon} \rvert^{p} \bigr] = p \lvert (\nabla v)_{\varepsilon} \rvert^{p - 2} (\nabla v)_{\varepsilon} \cdot \partial_{i} (\nabla v)_{\varepsilon} \longrightarrow p \lvert \nabla v \rvert^{p - 2} \nabla v \cdot \nabla(\partial_{i} v) \qquad \text{in } L^{k}_{\mathrm{loc}}(U).
\]
Hence $\lvert \nabla v \rvert^{p} \in W^{1, k}_{\mathrm{loc}}(U)$ with $\partial_{i} \lvert \nabla v \rvert^{p} = p \lvert \nabla v \rvert^{p - 2} \nabla v \cdot \nabla(\partial_{i} v)$. By \eqref{e:w-equation}, $\Delta v = \varsigma^{-1} (\lvert \nabla v \rvert^{p} + \lambda v - F) \in W^{1, k}_{\mathrm{loc}}(U)$. Thus $v_{i} := \partial_{i} v \in W^{1, k}_{\mathrm{loc}}(U)$ satisfies $\Delta v_{i} = \partial_{i} \Delta v \in L^{k}_{\mathrm{loc}}(U)$ in the sense of distributions. The same local Calder\'{o}n--Zygmund regularity argument as in Step 1 of Lemma~\ref{l:critical-blowup}, combined with mollification, gives $v_{i} \in W^{2, k}_{\mathrm{loc}}(U)$. Differentiating \eqref{e:w-equation} in the sense of distributions gives \eqref{e:w-differentiated}; since all its terms belong to $L^{k}_{\mathrm{loc}}(U)$, it holds almost everywhere.
\end{proof}

The following lemma establishes local uniform convergence of gradients without requiring convergence of $\nabla F_{j}$. Its proof combines the convexity of the Hamiltonian with estimates along trajectories driven by the formal minimizing feedback. For other approaches to vanishing-viscosity limits of Hamilton--Jacobi equations, see \cite{DNT26, lE10, hT11}.
\begin{lem} \label{l:gradient-convergence}
Let $U, U_{j} \subset \mathbb{R}^{N}$ be open, and suppose that every compact subset of $U$ is contained in $U_{j}$ for all sufficiently large $j$. Let $\varsigma_{j} > 0$ and $\lambda_{j} \ge 0$ with $\varsigma_{j} \to 0$ and $\lambda_{j} \to 0$. Let $F_{j} \in W^{1, \infty}_{\mathrm{loc}}(U_{j})$, and let $v_{j} \in W^{2, k}_{\mathrm{loc}}(U_{j})$ for all $k \in (1, \infty)$ satisfy
\begin{equation} \label{e:vj}
-\varsigma_{j} \Delta v_{j} + \lvert \nabla v_{j} \rvert^{p} + \lambda_{j} v_{j} = F_{j} \qquad \text{a.e. in } U_{j}.
\end{equation}
Let $v \in C^{2}(U)$ and $F \in C(U)$ satisfy $\lvert \nabla v \rvert^{p} = F$ in $U$, and assume that, for every compact $K \subset U$,
\begin{equation*}
\lVert v_{j} - v \rVert_{L^{\infty}(K)} + \lVert F_{j} - F \rVert_{L^{\infty}(K)} \stackrel{{j \to \infty}}{\longrightarrow} 0, \qquad \limsup_{j \to \infty} \bigl(\lVert \nabla v_{j} \rVert_{L^{\infty}(K)} + \lVert \nabla F_{j} \rVert_{L^{\infty}(K)} \bigr) < \infty.
\end{equation*}
Then $\nabla v_{j} \to \nabla v$ locally uniformly in $U$.
\end{lem}

\begin{proof}
{\it Step 1: localization and the feedback diffusion}. Fix a compact set $K \subset U$, and choose $\varrho_{1} \in (0, 1]$ such that $K_{1} := \{y : \operatorname{dist}(y, K) \le \varrho_{1}\} \subset U$. After discarding finitely many indices, $K_{1} \subset U_{j}$ and
\[
M_{\nabla} := \lVert \nabla v \rVert_{L^{\infty}(K_{1})} + \sup_{j} \, \lVert \nabla v_{j} \rVert_{L^{\infty}(K_{1})} < \infty, \qquad M_{F} := \sup_{j} \, \lVert \nabla F_{j} \rVert_{L^{\infty}(K_{1})} < \infty.
\]
Set
\[
\Lambda_{v} := \sup_{K_{1}} \, \lVert D^{2} v \rVert_{\mathrm{op}}, \qquad M_{b} := p M_{\nabla}^{p - 1}.
\]
By Lemma~\ref{l:w3s}, $v_{j} \in W^{3, k}_{\mathrm{loc}}(U_{j})$ for all $k \in (1, \infty)$. In particular, each $\partial_{i} v_{j}$ belongs to $W^{2, k}_{\mathrm{loc}}(U_{j})$ for all $k \in (1, \infty)$, and $\nabla v_{j}$ is continuous. Define the formal minimizing feedback and the gradient error
\[
b_{j} := -p \lvert \nabla v_{j} \rvert^{p - 2} \nabla v_{j}, \qquad \Psi_{j} := \nabla v_{j} - \nabla v \qquad \text{on } K_{1},
\]
and note that $\lvert b_{j} \rvert \le M_{b}$ on $K_{1}$. By $L(b_{j}) = (p - 1) \lvert \nabla v_{j} \rvert^{p}$, \eqref{e:w-differentiated}, and \eqref{e:vj}, almost everywhere on $K_{1}$,
\begin{align}
\varsigma_{j} \Delta v_{j} + b_{j} \cdot \nabla v_{j} - \lambda_{j} v_{j} & = -L(b_{j}) - F_{j}, \label{e:gc-generator} \\
\varsigma_{j} \Delta(\partial_{i} v_{j}) + b_{j} \cdot \nabla(\partial_{i} v_{j}) - \lambda_{j} \partial_{i} v_{j} & = -\partial_{i} F_{j}, \qquad i = 1, \dots, N. \label{e:gc-costate}
\end{align}

Fix $j$ and $x \in K$, and choose $T \in (0, 1]$ so small that $M_{b} T \le \varrho_{1} / 2$. The drift $\widetilde{b}_{j} := b_{j} \mathbf{1}_{K_{1}}$ is bounded and Borel, so there exists a weak solution of
\[
X_{t} = x + \int_{0}^{t} \widetilde{b}_{j}(X_{s}) \, \mathrm{d}s + \sqrt{2 \varsigma_{j}} \, W_{t}, \qquad t \in [0, T],
\]
where $W$ is an $N$-dimensional Brownian motion, see \cite[Prop.~5.3.6]{KS91}. Set $\tau := T \wedge \inf\{t \ge 0 : \lvert X_{t} - x \rvert \ge \varrho_{1}\}$. On $[0, \tau]$ the process stays in $\overline{B}(x, \varrho_{1}) \subset K_{1} \subset U_{j}$ and has drift $b_{j}(X_{t})$. By Lemma~\ref{l:w3s}, the functions $v_{j}$ and $\partial_{i} v_{j}$ belong to $W^{2, k}(B(x, \varrho_{1}))$ for all $k \in (1, \infty)$, and their gradients are continuous on $\overline{B}(x, \varrho_{1})$. Since the drift is bounded and $\varsigma_{j} > 0$, the It\^{o}--Krylov formula \cite[Th.~2.10.1]{nK80} applies to $v_{j}$ and $\partial_{i} v_{j}$ up to $\tau$ and justifies the use of the a.e.\ identities \eqref{e:gc-generator} and \eqref{e:gc-costate}. The resulting stochastic integrals have zero expectation because their integrands are bounded. Hence, for $t \in [0, T]$,
\begin{align} \label{e:gc-value}
\mathbb{E} \bigl[e^{-\lambda_{j} (t \wedge \tau)} v_{j}(X_{t \wedge \tau}) \bigr] & = v_{j}(x) - \mathbb{E} \int_{0}^{t \wedge \tau} e^{-\lambda_{j} s} \bigl(L(b_{j}) + F_{j} \bigr)(X_{s}) \, \mathrm{d}s, 
\\ \label{e:gc-costate-identity}
\mathbb{E} \bigl[e^{-\lambda_{j} (t \wedge \tau)} \partial_{i} v_{j}(X_{t \wedge \tau}) \bigr] & = \partial_{i} v_{j}(x) - \mathbb{E} \int_{0}^{t \wedge \tau} e^{-\lambda_{j} s} \, \partial_{i} F_{j}(X_{s}) \, \mathrm{d}s, \qquad i = 1, \dots, N.
\end{align}
Since $v \in C^{2}(U)$, the classical It\^{o}'s formula gives, for $t \in [0, T]$,
\begin{equation} \label{e:gc-ito-v}
\mathbb{E} \bigl[e^{-\lambda_{j} (t \wedge \tau)} v(X_{t \wedge \tau}) \bigr] = v(x) + \mathbb{E} \int_{0}^{t \wedge \tau} e^{-\lambda_{j} s} \bigl(\varsigma_{j} \Delta v + b_{j} \cdot \nabla v - \lambda_{j} v \bigr)(X_{s}) \, \mathrm{d}s.
\end{equation}

{\it Step 2: convexity estimate}. On $K_{1}$ define
\[
G_{j} := \lvert \nabla v \rvert^{p} - \lvert \nabla v_{j} \rvert^{p} + p \lvert \nabla v_{j} \rvert^{p - 2} \nabla v_{j} \cdot \Psi_{j} = L(b_{j}) + b_{j} \cdot \nabla v + \lvert \nabla v \rvert^{p}.
\]
Integrating $\frac{\mathrm{d}}{\mathrm{d}\zeta} \lvert \nabla v_{j} - \zeta \Psi_{j} \rvert^{p}$ over $\zeta \in [0, 1]$ gives
\[
G_{j} = p \int_{0}^{1} \bigl\langle \lvert \nabla v_{j} - \zeta \Psi_{j} \rvert^{p - 2} (\nabla v_{j} - \zeta \Psi_{j}) - \lvert \nabla v_{j} \rvert^{p - 2} \nabla v_{j}, -\Psi_{j} \bigr\rangle \, \mathrm{d}\zeta.
\]
For $\zeta \in (0, 1]$, the inequality \cite[Sect.~10, (I)]{pL06}, applied to the vectors $\nabla v_{j} - \zeta \Psi_{j}$ and $\nabla v_{j}$, bounds the integrand from below by $2^{2 - p} \zeta^{p - 1} \lvert \Psi_{j} \rvert^{p}$. Hence
\begin{equation} \label{e:bregman}
G_{j} \ge c_{p} \lvert \Psi_{j} \rvert^{p} \quad \text{on } K_{1}, \qquad c_{p} := 2^{2 - p}.
\end{equation}
Subtracting \eqref{e:gc-ito-v} from \eqref{e:gc-value} with $t = T$, so that $t \wedge \tau = \tau$, and using $F = \lvert \nabla v \rvert^{p}$, we obtain
\[
(v_{j} - v)(x) - \mathbb{E} \bigl[e^{-\lambda_{j} \tau} (v_{j} - v)(X_{\tau}) \bigr] = \mathbb{E} \int_{0}^{\tau} e^{-\lambda_{j} s} \bigl(G_{j} + F_{j} - F + \varsigma_{j} \Delta v - \lambda_{j} v \bigr)(X_{s}) \, \mathrm{d}s.
\]
Set
\begin{equation*}
\eta_{j} := 2 \lVert v_{j} - v \rVert_{L^{\infty}(K_{1})} + \lVert F_{j} - F \rVert_{L^{\infty}(K_{1})} + \varsigma_{j} \lVert \Delta v \rVert_{L^{\infty}(K_{1})} + \lambda_{j} \lVert v \rVert_{L^{\infty}(K_{1})}.
\end{equation*}
By our assumptions, $\eta_{j} \to 0$ as $j \to \infty$. Since $T \le 1$, we have $\mathbb{E} \int_{0}^{\tau} e^{-\lambda_{j} s} G_{j}(X_{s}) \, \mathrm{d}s \le \eta_{j}$, and therefore, by \eqref{e:bregman},
\begin{equation} \label{e:bregman-traj}
\mathbb{E} \int_{0}^{\tau} \lvert \Psi_{j}(X_{t}) \rvert^{p} \, \mathrm{d}t \le c_{p}^{-1} e^{\lambda_{j} T} \eta_{j}.
\end{equation}

{\it Step 3: gradient comparison}. Fix $t \in [0, T]$. Since $X_{s} \in K_{1}$ for $s \le \tau$, $\lvert \partial_{i} F_{j} \rvert \le M_{F}$ and $\lvert \partial_{i} v_{j} \rvert \le M_{\nabla}$ on $K_{1}$, and $0 \le 1 - e^{-\lambda_{j} (t \wedge \tau)} \le \lambda_{j} t$, the identity \eqref{e:gc-costate-identity} gives
\[
\bigl\lvert \partial_{i} v_{j}(x) - \mathbb{E} \, \partial_{i} v_{j}(X_{t \wedge \tau}) \bigr\rvert \le (M_{F} + M_{\nabla} \lambda_{j}) t, \qquad i = 1, \dots, N. 
\]
Combining the coordinates,
\begin{equation} \label{e:costate-bound}
\bigl\lvert \nabla v_{j}(x) - \mathbb{E} \, \nabla v_{j}(X_{t \wedge \tau}) \bigr\rvert \le \sqrt{N} (M_{F} + M_{\nabla} \lambda_{j}) t.
\end{equation}

Since $X_{t \wedge \tau} = X_{\tau}$ for $t \ge \tau$, we have pathwise 
\[
\frac{1}{T} \int_{0}^{T} \lvert \Psi_{j}(X_{t \wedge \tau}) \rvert \, \mathrm{d}t = \frac{1}{T} \int_{0}^{\tau} \lvert \Psi_{j}(X_{t}) \rvert \, \mathrm{d}t + \frac{T - \tau}{T} \, \lvert \Psi_{j}(X_{\tau}) \rvert.
\]
The last term is at most $M_{\nabla} \mathbf{1}_{\{\tau < T\}}$, since $X_{\tau} \in K_{1}$. For the first term, H\"{o}lder's inequality and \eqref{e:bregman-traj} give
\[
\frac{1}{T} \, \mathbb{E} \int_{0}^{\tau} \lvert \Psi_{j}(X_{t}) \rvert \, \mathrm{d}t \le \frac{1}{T} \Bigl(\mathbb{E} \int_{0}^{\tau} \lvert \Psi_{j}(X_{t}) \rvert^{p} \, \mathrm{d}t \Bigr)^{1 / p} T^{1 - 1 / p} \le T^{-1 / p} \bigl(c_{p}^{-1} e^{\lambda_{j} T} \eta_{j} \bigr)^{1 / p}.
\]
Taking expectations in the pathwise identity, we obtain
\begin{equation} \label{e:gc-Psi}
\frac{1}{T} \int_{0}^{T} \mathbb{E} \lvert \Psi_{j}(X_{t \wedge \tau}) \rvert \, \mathrm{d}t \le T^{-1 / p} \bigl(c_{p}^{-1} e^{\lambda_{j} T} \eta_{j} \bigr)^{1 / p} + M_{\nabla} \, \mathbb{P}(\tau < T).
\end{equation}

Next, let $t \in [0, T]$. The segment from $x$ to $X_{t \wedge \tau}$ lies in the convex set $\overline{B}(x, \varrho_{1}) \subset K_{1}$, so $\lvert \nabla v(X_{t \wedge \tau}) - \nabla v(x) \rvert \le \Lambda_{v} \lvert X_{t \wedge \tau} - x \rvert$. Since $\lvert b_{j} \rvert \le M_{b}$ on $K_{1}$ and $\mathbb{E} \lvert W_{t \wedge \tau} \rvert^{2} = N \, \mathbb{E}(t \wedge \tau) \le N T$ by optional stopping,
\begin{equation} \label{e:gc-displacement}
\mathbb{E} \, \lvert \nabla v(X_{t \wedge \tau}) - \nabla v(x) \rvert \le \Lambda_{v} \Bigl(M_{b} T + \sqrt{2 \varsigma_{j}} \, \bigl(\mathbb{E} \lvert W_{t \wedge \tau} \rvert^{2} \bigr)^{1 / 2} \Bigr) \le \Lambda_{v} \bigl(M_{b} T + \sqrt{2 N \varsigma_{j} T} \bigr).
\end{equation}

On $\{\tau < T\}$ we have $\lvert X_{\tau} - x \rvert = \varrho_{1}$ and $M_{b} T \le \varrho_{1} / 2$, so $\sup_{t \le T} \lvert W_{t} \rvert \ge \varrho_{1} / (2 \sqrt{2 \varsigma_{j}})$. By Chebyshev's and Doob's inequalities,
\begin{equation} \label{e:gc-exit}
\mathbb{P}(\tau < T) \le \frac{8 \varsigma_{j}}{\varrho_{1}^{2}} \, \mathbb{E} \bigl( \sup_{t \le T} \lvert W_{t} \rvert \bigr) ^{2} \le \frac{32 N \varsigma_{j} T}{\varrho_{1}^{2}}.
\end{equation}

{\it Step 4: conclusion}. For $t \in [0, T]$,
\[
\lvert \nabla v_{j}(x) - \nabla v(x) \rvert \le \bigl\lvert \nabla v_{j}(x) - \mathbb{E} \, \nabla v_{j}(X_{t \wedge \tau}) \bigr\rvert + \mathbb{E} \lvert \Psi_{j}(X_{t \wedge \tau}) \rvert + \mathbb{E} \lvert \nabla v(X_{t \wedge \tau}) - \nabla v(x) \rvert.
\]
Averaging over $t \in [0, T]$ and using \eqref{e:costate-bound}, \eqref{e:gc-Psi}, \eqref{e:gc-displacement} and \eqref{e:gc-exit}, we obtain, for every $x \in K$ and for all sufficiently small $T > 0$,
\[
\begin{aligned}
	\lvert \nabla v_{j}(x) - \nabla v(x) \rvert \le & \sqrt{N} \, \frac{(M_{F} + M_{\nabla} \lambda_{j}) T}{2} + T^{-1 / p} \bigl(c_{p}^{-1} e^{\lambda_{j} T} \eta_{j} \bigr)^{1 / p} \\
							&+ \frac{32 M_{\nabla} N \varsigma_{j} T}{\varrho_{1}^{2}} + \Lambda_{v} \bigl(M_{b} T + \sqrt{2 N \varsigma_{j} T} \bigr).
\end{aligned}
\] 
Since the estimate is uniform in $x \in K$, for all sufficiently small $T > 0$ we have $\limsup_{j \to \infty} \lVert \nabla v_{j} - \nabla v \rVert_{L^{\infty}(K)} \le (\frac{1}{2} \sqrt{N} M_{F} + \Lambda_{v} M_{b}) T$, and letting $T \downarrow 0$ concludes the proof.
\end{proof}

We can now complete the analysis of the supercritical blow-up limits.
\begin{cor} \label{c:super-C1}
Assume Condition~\ref{cond:supercritical}. Then $u_{j} \to \Phi$ in $C^{1}_{\mathrm{loc}}(H)$.
\end{cor}
\begin{proof}
By Lemmas~\ref{l:blowup} and~\ref{l:super-profile}, $\Phi$ is given by \eqref{e:super-unique-blowup}; in particular, $\Phi \in C^{\infty}(H)$, $\lvert \nabla \Phi \rvert^{p} = a_{0} \rho^{-\beta}$, and $u_{j} \to \Phi$ locally uniformly in $H$ along the original sequence. We apply Lemma~\ref{l:gradient-convergence} with
\[
U = H, \quad U_{j} = \Omega_{j}, \quad v_{j} = u_{j}, \quad v = \Phi, \quad F = a_{0} \rho^{-\beta}, \quad \varsigma_{j} = r_{j}^{m - 2}, \quad \lambda_{j} = \lambda r_{j}^{m}.
\]
Its hypotheses follow from Lemma~\ref{l:blowup}, \eqref{e:LL-regularity}, \eqref{e:scaled-PDE}, \eqref{e:blowup-limit}, \eqref{e:forcing-limit} and \eqref{e:scaled-gradient}. For $\nabla F_{j}$, note that $\nabla F_{j}(y) = r_{j}^{\beta + 1} \nabla f(r_{j} y)$, so \eqref{e:wdf} gives $F_{j} \in W^{1, \infty}_{\mathrm{loc}}(\Omega_{j})$ and $\lvert \nabla F_{j} \rvert \le C \, d_{j}^{-\beta - 1}$, which is bounded on compact subsets of $H$. Hence $\nabla u_{j} \to \nabla \Phi$ locally uniformly in $H$, completing the proof.
\end{proof}

\section{Conclusion of the proofs} \label{s:cop}

\subsection{Proof of Theorems~\ref{t:c}, \ref{t:supc} and Corollary~\ref{c:main}}
Set $\mathcal{A}(s) := A(s)$ under Condition~\ref{cond:critical} and $\mathcal{A}(s) := s^{1 / p}$ under Condition~\ref{cond:supercritical}, for $s \ge 0$. For each $x_{0} \in \partial \Omega$, choose a rotation $Q_{x_{0}}$ such that $Q_{x_{0}} e_{N} = -n(x_{0})$, and define
\[
u_{x_{0}, r}(\xi) := r^{-\alpha} \bigl(u(x_{0} + r Q_{x_{0}} \xi) - u(x_{0}) \bigr), \qquad \xi \in \Omega_{x_{0}, r},
\]
where $\Omega_{x_{0}, r} := r^{-1} Q_{x_{0}}^{\top}(\Omega - x_{0})$. By the $C^{2}$ regularity and compactness of $\partial \Omega$, every compact $K \Subset H$ is contained in $\Omega_{x_{0}, r}$ for all sufficiently small $r$, uniformly in $x_{0}$.

We now combine the preceding convergence results uniformly over $x_{0} \in \partial \Omega$. Theorems~\ref{t:c}, \ref{t:supc} and Corollary~\ref{c:main} follow immediately from the following proposition by evaluating at $e_{N}$, using $x = \pi_{x} - d(x)n(\pi_{x})$, and recalling that $A(0) = (q - 1)^{q - 1}$.

\begin{prop}
For every $K \Subset H$,
\begin{equation} \label{e:uniform-blowup}
\lim_{r \downarrow 0} \sup_{x_{0} \in \partial \Omega} \lVert u_{x_{0}, r} - \Phi_{x_{0}} \rVert_{C^{1}(K)} = 0,
\end{equation}
where, for $(z, \rho) \in H$,
\begin{equation} \label{e:unified-blowup}
\Phi_{x_{0}}(z, \rho) := -\frac{\mathcal{A}(a(x_{0}))}{\alpha} \rho^{\alpha}, \qquad \nabla \Phi_{x_{0}}(z, \rho) = -\mathcal{A}(a(x_{0})) \rho^{-\theta} e_{N}.
\end{equation}
\end{prop}
\begin{proof}
The constants $C_{H}$ and $C_{G}$ in Proposition~\ref{p:LL} are independent of the boundary point $x_{0}$. By the $C^{2}$ regularity and compactness of $\partial \Omega$, the graph neighborhood $\mathcal{C}_{0}$ and the constant $C_{g}$ in \eqref{e:g} can be chosen uniformly in $x_{0}$. Consequently, the graph representation and the estimates \eqref{e:scaled-graph} in Lemma~\ref{l:blowup} hold uniformly in $x_{0}$. The constant $\delta_{0}$ and the bound $C_{d}$ on $\lvert \Delta d \rvert$ are also independent of $x_{0}$.

The constants $C_{E}$ and $C_{X}$ in Lemma~\ref{l:Et} are independent of $x_{0}$, and, for fixed $\xi$ and $T$, the remainder terms in \eqref{e:Ea} are $O(\varepsilon_{j} + r_{j}^{\beta})$ uniformly in $x_{0}$, by \eqref{e:scaled-growth} and the boundedness of $u$. Moreover, the proof of Lemma~\ref{l:scaled-occupation} gives, for fixed $\xi$, $T$ and $M$, an occupation bound independent of $x_{0}$, since the constants $c_{0}$ and $C_{0}$ in the critical case and the positive lower bound $a_{*}$ in the supercritical case are independent of $x_{0}$.

For fixed $\xi$, $h$, $T$ and $M$, the compact set $\mathcal{K}$ and the cylinder $\mathcal{C}$ in the proof of Proposition~\ref{p:ti} can be chosen independently of $x_{0}$. By the uniform graph estimates, $\ell_{j} = O(r_{j})$ uniformly in $x_{0}$, so that $\delta_{j} \to 0$ and $\gamma_{j} \to 0$ uniformly in $x_{0}$. Finally,
\begin{equation} \label{e:a}
\lim_{\varrho \downarrow 0} \sup_{x_{0} \in \partial \Omega} \operatorname*{ess\,sup}_{x \in \Omega \cap B(x_{0}, \varrho)} \lvert a(x) - a(x_{0}) \rvert = 0.
\end{equation}
This follows from \eqref{e:critical-trace} and the uniform continuity of the boundary trace in the critical case, and from the uniform continuity of $a$ in the supercritical case. In particular, the errors $\omega_{j}$ in the proof of Proposition~\ref{p:ti} are bounded by the supremum in \eqref{e:a} with $\varrho = r_{j} \bigl(1 + \sup_{y \in \mathcal{K}} \lvert y \rvert \bigr)$, and hence tend to zero uniformly in $x_{0}$. Together with the boundedness of $a$ on $\partial \Omega$, the preceding observations show that all the estimates in the tangential invariance argument remain valid for sequences of moving boundary points.

Consider arbitrary sequences $x_{0}^{j} \in \partial \Omega$ and $r_{j} \downarrow 0$, and write $u_{j} := u_{x_{0}^{j}, r_{j}}$ and $\Omega_{j} := \Omega_{x_{0}^{j}, r_{j}}$. For these blow-ups, the rescaled distance and state cost are
\[
d_{j}(y) := \frac{d(x_{0}^{j} + r_{j} Q_{x_{0}^{j}} y)}{r_{j}}, \qquad F_{j}(y) := r_{j}^{\beta} \bigl(f(x_{0}^{j} + r_{j} Q_{x_{0}^{j}} y) - \lambda u(x_{0}^{j}) \bigr).
\]
After passing to a subsequence, $x_{0}^{j} \to x_{0}^{*} \in \partial \Omega$. By the uniform H\"{o}lder and boundary graph estimates, as in Lemma~\ref{l:blowup}, we may assume that $u_{j} \to \Phi$ locally uniformly in $H$. By \eqref{e:scaled-holder} and \eqref{e:scaled-growth}, $\Phi$ extends continuously to $\overline{H}$ with $\Phi(0) = 0$. By the preceding observations, the proof of Proposition~\ref{p:ti} applies to these blow-ups, with $a(x_{0}^{j})$ in place of $a_{0}$, and yields the tangential invariance of $\Phi$. Moreover, the uniform boundary growth estimate \eqref{e:asn} gives \eqref{e:pa} for $\Phi$.

The uniform graph estimates give $d_{j}(z, \rho) \to \rho$ locally uniformly in $H$. Since $a(x_{0}^{j}) \to a(x_{0}^{*})$, \eqref{e:a} and the boundedness of $u$ yield
\[
F_{j}(z, \rho) \longrightarrow a(x_{0}^{*}) \rho^{-\beta} \qquad \text{in } L^{\infty}_{\mathrm{loc}}(H).
\]
In the supercritical case, this convergence is locally uniform.

Under Condition~\ref{cond:critical}, the local $W^{2, k}$ estimates in the proof of Lemma~\ref{l:critical-blowup}, for $k > N$, follow from \eqref{e:scaled-growth}, \eqref{e:scaled-gradient} and the rescaled equation \eqref{e:scaled-PDE}, using the local uniform bounds on $F_{j}$. Hence $u_{j} \to \Phi$ in $C^{1}_{\mathrm{loc}}(H)$, and the classification argument in that lemma gives $\Phi = \Phi_{x_{0}^{*}}$.

Under Condition~\ref{cond:supercritical}, the proof of Lemma~\ref{l:super-profile} gives $\Phi = \Phi_{x_{0}^{*}}$. Moreover, the bounds \eqref{e:scaled-gradient} and \eqref{e:wdf} give local uniform bounds for $\nabla u_{j}$ and $\nabla F_{j}$ in $H$, independently of $x_{0}^{j}$. Indeed, on every compact subset of $H$, for all sufficiently large $j$,
\[
\lvert \nabla u_{j} \rvert \le C_{G} d_{j}^{-\theta}, \qquad \lvert \nabla F_{j} \rvert \le C d_{j}^{-\beta - 1} \quad \text{a.e.},
\]
where $C$ is the bound in \eqref{e:wdf}. We may therefore apply Lemma~\ref{l:gradient-convergence} as in the proof of Corollary~\ref{c:super-C1}, obtaining $u_{j} \to \Phi_{x_{0}^{*}}$ in $C^{1}_{\mathrm{loc}}(H)$.

Since $\mathcal{A} \circ a$ is continuous on $\partial \Omega$, \eqref{e:unified-blowup} gives $\Phi_{x_{0}^{j}} \to \Phi_{x_{0}^{*}}$ in $C^{1}_{\mathrm{loc}}(H)$. Consequently, along the extracted subsequence, $\lVert u_{x_{0}^{j}, r_{j}} - \Phi_{x_{0}^{j}} \rVert_{C^{1}(K)} \to 0$ for every $K \Subset H$. Since the original sequences were arbitrary, this proves \eqref{e:uniform-blowup}.
\end{proof}

\subsection{Proof of Theorems \ref{t:critical-feedback} and \ref{t:supercritical-feedback}}
As observed by Lasry and Lions \cite[Sect.~VII.2 and VII.4]{LL89}, the sharp boundary gradient asymptotics yield the admissibility and optimality of the minimizing feedbacks by standard arguments. Here the asserted asymptotics of $b^{*}$ follow directly from \eqref{e:cmg}, \eqref{e:smg} and \eqref{def:b0}. We briefly verify the remaining claims.

{\it Admissibility}. 
These asymptotics imply that, for some $\delta_{1} \in (0, \delta_{0})$,
\begin{equation} \label{e:drift-inward}
	d(x) b^{*}(x) \cdot \nabla d(x) \ge \frac{1 + q}{2}, \qquad x \in \Omega_{\delta_{1}}.
\end{equation}
Choose $\psi \in C^{2}(\Omega)$ with $\psi = -\log d$ on $\Omega_{\delta_{1}}$. By \eqref{e:drift-inward}, for some $C_{d}, C > 0$,
\[
\Delta \psi + b^{*} \cdot \nabla \psi
\le \frac{1 - q}{2d^{2}} + \frac{C_{d}}{d} \le C
\qquad \text{on } \Omega_{\delta_{1}}.
\]
By continuity on the compact complement, this bound holds throughout $\Omega$, where $\psi \ge -C_{1}$ for some $C_{1} > 0$. Localization and \cite[Prop.~5.3.6, Cor.~5.3.11]{KS91} yield a weak solution $X$ of \eqref{e:feedback-SDE}, unique in law up to the exit time. For large $i$, set $\tau_{i} := \inf\{t \ge 0 : d(X_{t}) \le i^{-1}\}$. It\^{o}'s formula gives
\[
\mathbb{P}(\tau_{i} \le t) \log i \le \psi(x) + C t + C_{1}.
\]
Letting $i \to \infty$ shows that $X$ remains in $\Omega$ almost surely. Thus the feedback $b^{*}$ is admissible.

{\it Optimality}. Since $b^{*}$ realizes equality in \eqref{e:Fenchel-duality} with $v = \nabla u$, the proof of Proposition~\ref{p:verification} gives equality in \eqref{e:uu} for $b_{t} = b^{*}(X_{t})$ and $\tau = T$. Letting $T \to \infty$ and using the boundedness of $u$ and monotone convergence for the running cost after adding $c_{f}$, we obtain $J_{\lambda}(x; b^{*}(X)) = u(x) = V_{\lambda}(x)$. Thus $b^{*}$ is optimal.

\subsection{Proof of Corollary~\ref{c:ergodic}}
Fix $\lambda_{0} > 0$ and set $f_{0} := f - \kappa + \lambda_{0} u_{0}$.
Then
\[
-\Delta u_{0} + \lvert \nabla u_{0} \rvert^{p} + \lambda_{0} u_{0} = f_{0}
\qquad \text{in } \Omega.
\]
By \cite[Proof of Th.~VI.3]{LL89}, $u_{0}$ is the maximal solution of this discounted equation. Moreover, the boundedness of $u_{0}$ ensures that $f_{0}$ is bounded from below and has the same boundary coefficient as $f$. In the supercritical case, the gradient estimate \cite[Th.~IV.1]{LL89} also ensures that $f_{0}$ satisfies \eqref{e:wdf}. Thus $f_{0}$ satisfies the corresponding condition. The boundary asymptotics follow from Theorem~\ref{t:c} or~\ref{t:supc}, and the preceding admissibility argument applies to $b_{0}^{*}$. Thus, for every $x \in \Omega$, the feedback equation with $b_{0}^{*}$ in place of $b^{*}$ admits a global weak solution $X$, unique in law, and the control $t \mapsto b_{0}^{*}(X_{t})$ belongs to $\operatorname{Adm}(x)$.

Finally, using \eqref{e:Fenchel-duality} and \eqref{e:ergodic}, the same verification argument as in the proof of Proposition~\ref{p:verification}, now without discounting, gives, for every $b \in \operatorname{Adm}(x)$ and $T > 0$,
\[
\frac{1}{T} \mathbb{E} \int_{0}^{T} \bigl(L(b_{t}) + f(X_{t}) \bigr) \, \mathrm{d}t \ge \kappa + \frac{u_{0}(x) - \mathbb{E} \, u_{0}(X_{T})}{T},
\]
with equality for $b_{t} = b_{0}^{*}(X_{t})$. Letting $T \to \infty$ and using the boundedness of $u_{0}$, we obtain \eqref{e:ergodic-value} and the optimality of $b_{0}^{*}$.

\end{document}